\documentclass[12pt]{amsart}
\usepackage{amssymb,latexsym,pstricks}

\numberwithin{equation}{section}

\newtheorem{theorem}{Theorem}[section]
\newtheorem{proposition}[theorem]{Proposition}

\newtheorem{corollary}[theorem]{Corollary}
\newtheorem{lemma}[theorem]{Lemma}

\theoremstyle{definition}

\newtheorem{remark}[theorem]{Remark}
\newtheorem{example}[theorem]{Example}
\newtheorem{definition}[theorem]{Definition}

\let\oldmarginpar\marginpar
\renewcommand\marginpar[1]{\-\oldmarginpar[\raggedleft\small\sf
#1]{\raggedright\small\sf #1}}

\newcommand{\myAA}{\mathcal{A}}

\newcommand{\BB}{\mathcal{B}}

\newcommand{\PP}{\mathbb{P}}
\newcommand{\FF}{\mathcal{F}}

\newcommand{\ZZ}{\mathbb{Z}}

\newcommand{\RR}{\mathbb{R}}
\newcommand{\QQ}{\mathbb{Q}}

\newcommand{\mat}[4]{\left(\!\!\begin{array}{cc}
#1 & #2 \\ #3 & #4 \\
\end{array}\!\!\right)}
\newcommand{\rem}[1]{\left\langle#1\right\rangle}

\begin{document}

\title[Positivity in rank 2 cluster algebras]
{Positivity and canonical bases in rank $2$ cluster algebras of finite and
affine types}
%Canonical Bases of Finite and Affine Cluster Algebras of Rank 2}
\author{Paul Sherman}
\address{\noindent Department of Mathematics, Northeastern University,
Boston, MA 02115, USA}
\email{paul@avalon.math.neu.edu}

\author{Andrei Zelevinsky}
\address{\noindent Department of Mathematics, Northeastern University,
Boston, MA 02115, USA}
\email{andrei@neu.edu}

\begin{abstract}
The main motivation for the study of cluster algebras initiated
in~\cite{fz-ClusterI,fz-ClusterII,fz-ClusterIII}
%math.RT/0104151, math.RA/0208229 and math.RT/0305434
was to design an
algebraic framework for understanding total positivity
and canonical bases in semisimple algebraic groups.
In this paper, we introduce and explicitly construct the canonical basis for
a special family of cluster algebras of rank~$2$.
\end{abstract}

\date{July 7, 2003; revised December 12. 2003}

\dedicatory{To Borya Feigin on the occasion of his jubilee}

\thanks{Research %The authors were
supported in part by a US Department of Education GAANN grant (P.S.)
and NSF (DMS) grant \# 0200299 (A.Z.).}

\maketitle

\medskip

\begin{flushright}
\begin{minipage}{4.30in}
\small\textsl{
\noindent
{\bf ju-bi-lee}
\hfill
\hfill
\linebreak
\smallskip
\noindent
1 : a year of emancipation and restoration provided by ancient Hebrew law
to be kept every 50 years by the emancipation of Hebrew slaves, restoration
of alienated lands to their former owners, and omission of all cultivation
of the land
\linebreak
\smallskip
\noindent
2 a : a special anniversary; especially : a 50th anniversary
\linebreak
\noindent
b : a celebration of such an anniversary
\begin{flushright}
{\rm Merriam-Webster Dictionary}
\end{flushright}
}
\end{minipage}
\end{flushright}

\medskip

\tableofcontents

\section{Introduction}

This paper continues the study of cluster algebras initiated
in~\cite{fz-ClusterI,fz-ClusterII,fz-ClusterIII}. The principal
motivation for this theory was to design an algebraic framework
for understanding total positivity and canonical bases in
semisimple algebraic groups; see~\cite{zel-cambridge} for more
details. Since its inception, the theory of cluster algebras has
developed several interesting connections:

\begin{itemize}
\item Discrete dynamical systems~\cite{fz-Laurent}.

\item Al. Zamolodchikov's $Y$-systems in thermodynamic
Bethe Ansatz~\cite{yga}.

\item A new family of convex polytopes (generalized associahedra)
including as special cases Stasheff's associahedron, and
Bott-Taubes cyclohedron~\cite{yga,chap-fom-zel}.

\item Quiver representations~\cite{marsh-rei-zel}.

\item Grassmannians and projective configurations~\cite{scott}.

\item Poisson geometry and Teichmueller
spaces~\cite{fock-gonch1,fock-gonch2,gekh-shap-vai1,gekh-shap-vai2}.
\end{itemize}

Amid these exciting developments, virtually no progress has been
made (at least, documented) towards the original motivation.
This paper presents the first step in this direction
by introducing and explicitly constructing the canonical basis $\BB$
for a special class of rank~$2$ cluster algebras of finite anf
affine type.
More precisely, in Theorem~\ref{th:canonicalbasis} we show that~$\BB$
is uniquely determined by a certain positivity property.
Needless to say, we hope to be able one day to extend the results of this paper
to a much more general class of cluster algebras.

Restricting our attention to rank~$2$ cluster algebras allows us
to give all the necessary definitions in a quick and painless way.
Thus, the treatment in this paper is completely self-contained,
and no preliminary knowledge on cluster algebras is assumed from
the reader.

\section{Main results}

Let $\FF = \QQ(y_1,y_2)$ be the field of rational functions in two
(commuting) independent variables $y_1$ and $y_2$ with rational
coefficients.
Given positive integers $b$ and $c$, recursively define elements
$y_m \in \FF$ for $m \in \ZZ$ by the relations
\begin{equation}
\label{eq:clusterrelations}
y_{m-1} y_{m+1} =
\left\{
\begin{array}[h]{ll}
y_m^b + 1 & \quad \mbox{$m$ odd;} \\
y_m^c + 1 & \quad \mbox{$m$ even.}
\end{array}
\right.
\end{equation}
The main object of study in this paper is the \emph{cluster
algebra} $\myAA(b,c)$: by definition, this is a subring of $\FF$
generated by the $y_m$ for all $m \in \ZZ$.
The elements $y_m$ are called \emph{cluster variables} and the relations
(\ref{eq:clusterrelations}) are called the \emph{exchange relations}.
The sets $\{y_m,y_{m+1}\}$ for $m\in\ZZ$ are called \emph{clusters}.

In the terminology of~\cite{fz-ClusterI},  $\myAA(b,c)$ is a
(coefficient-free) cluster algebra of rank~$2$.
The most general definition of cluster algebras given in~\cite{fz-ClusterI}
allows terms in the right hand side of (\ref{eq:clusterrelations}) to have
nontrivial
coefficients.
In Section~\ref{sec:coeffs-removed}, we will show that for rank~$2$ cluster
algebras,
it suffices to restrict to the coefficient-free case since
the coefficients may be removed by a rescaling of the~$y_m$.

%The main motivation for the study of cluster algebras initiated
%in~\cite{fz-ClusterI,fz-ClusterII,fz-ClusterIII} was to design an
%algebraic framework for understanding \emph{total positivity}
%and \emph{canonical bases} in semisimple algebraic groups.
In this paper, we %make a first step in this direction by
introduce and explicitly construct the canonical basis $\BB$ in
$\myAA(b,c)$ for $bc \leq 4$.
More precisely, in Theorem~\ref{th:canonicalbasis} we show that~$\BB$
is uniquely determined by a certain positivity property.
We expect %Theorem~\ref{th:canonicalbasis}
this result to hold in much greater generality; in particular, in
a sequel to this paper we are planning to investigate whether the
assumption $bc \leq 4$ can be lifted.

Before stating the new results, we briefly discuss some known
properties of $\myAA(b,c)$.

\subsection{Laurent phenomenon}
It is clear from (\ref{eq:clusterrelations}) that
every cluster of $\myAA(b,c)$ is a transcendence basis of the
ambient field $\FF$, so every element of $\myAA(b,c)$ is uniquely
expressed as a rational function in $y_m$ and $y_{m+1}$, for every
$m \in \ZZ$.
According to the \emph{Laurent phenomenon} established in
\cite{fz-ClusterI,fz-Laurent}, all these rational functions are
actually Laurent polynomials with integer coefficients.
The following stronger result is a special case of the results
in~\cite{fz-ClusterIII}:
\begin{equation}
\label{eq:upper}
\myAA(b,c) = \bigcap_{m \in \ZZ} \ZZ[y_m^{\pm 1}, y_{m+1}^{\pm 1}]
= \bigcap_{m = 0}^2 \ZZ[y_m^{\pm 1}, y_{m+1}^{\pm 1}],
\end{equation}
where $\ZZ[y_m^{\pm 1}, y_{m+1}^{\pm 1}]$ denotes the ring of
Laurent polynomials with integer coefficients in $y_m$ and
$y_{m+1}$.  The symmetry of the exchange relations
(\ref{eq:clusterrelations}) allows the second intersection in
(\ref{eq:upper}) to be taken over any three consecutive
clusters.

\subsection{Standard monomial basis}
Another special case of the results in \cite{fz-ClusterIII} is the
following:
\begin{equation}
\label{eq:monomial-basis}
\text{the set $\{y_0^{a_0} y_1^{a_1} y_2^{a_2} y_3^{a_3} \ : \ a_m
\in \ZZ_{\geq 0}, a_0 a_2 = a_1 a_3 = 0\}$ is a $\ZZ$-basis of
$\myAA(b,c)$.}
\end{equation}
Again, in view of the obvious symmetry, the variables $y_0, \dots, y_3$
in (\ref{eq:monomial-basis}) can be replaced by any four
consecutive cluster variables.
We refer to the basis in (\ref{eq:monomial-basis}) as a
\emph{standard monomial basis}.

\subsection{Finite generation}
As an immediate corollary of (\ref{eq:monomial-basis})
(cf.~\cite{fz-ClusterIII}), we obtain that
the algebra $\myAA(b,c)$ is finitely generated and has the
following presentation:
\begin{equation}
\label{eq;presentation}
\myAA(b,c) = \ZZ[y_0, y_1, y_2, y_3]/\langle
y_0 y_2 - y_1^b -1, \, y_1 y_3 - y_2^c -1 \rangle \ .
\end{equation}

\subsection{Finite type classification}
The structure of $\myAA(b,c)$ is closely related
to the root system associated with the (generalized) Cartan matrix
\begin{equation}
\label{eq:cartan-rank2}
A=A(b,c)=\mat{2}{-b}{-c}{2} \ .
\end{equation}
Some properties of this root system will be reviewed in
Section~\ref{sec:general-rank2}.
Borrowing the terminology from the theory of Kac-Moody algebras
(see~\cite{kac}), we say that $\myAA(b,c)$ is of \emph{finite}
(resp.~\emph{affine, indefinite}) type if so is $A(b,c)$, i.e., if $bc \leq
3$
(resp. $bc = 4$, $bc > 4$).
The following result is proved in~\cite[Theorem~6.1]{fz-ClusterI}
(it is extended to cluster algebras of an arbitrary rank
in~\cite{fz-ClusterII}).

\begin{proposition}
\label{pr:finite-type-periodicity}
The algebra $\myAA(b,c)$ is of finite type if and only if
it has finitely many distinct cluster variables.
More precisely, if $bc \leq 3$ then
$$y_m = y_n \, \, \Leftrightarrow \,\,  m \equiv n \ \mod \ (h+2),$$
where $h$ is the Coxeter number given by Table~\ref{tab:Coxeter number}
below; whereas, if $bc \geq 4$ then all $y_m$ are distinct.
\end{proposition}

\begin{table}[h]
\begin{center}
\begin{tabular}{|c|c|c|c|c|}
\hline
$bc$  & 1 & 2 & 3 & $\geq 4$ \\
\hline
$h$  & 3 & 4 & 6 & $\infty$ \\
\hline
\end{tabular} \,.
\end{center}
\medskip
\caption{The Coxeter number}
\label{tab:Coxeter number}
\end{table}

\subsection{Positivity and canonical basis}
We now turn to the statements of the new results.
We will introduce the ``canonical" $\ZZ$-basis $\BB$ in~$\myAA(b,c)$,
in the cases when $\myAA(b,c)$ is of finite or affine type
(note that the basis in (\ref{eq:monomial-basis}) is not canonical,
because it involves an arbitrary choice of four consecutive
cluster variables).
Our construction of $\BB$ is based on the following concept of
\emph{positivity}.

\begin{definition}
\label{def:positivity}
A non-zero element $y \in \myAA(b,c)$ is \emph{positive} if
for every $m \in \ZZ$, all the coefficients
in the expansion of $y$ as a Laurent polynomial in $y_m$ and $y_{m+1}$ are
positive.
\end{definition}

The set of positive elements is a \emph{semiring}, i.e., is closed
under addition and multiplication.
The following theorem which is the main result of the paper gives
a complete description of this semiring in the finite and affine
type cases.

\begin{theorem}
\label{th:canonicalbasis}
Suppose that $bc \leq 4$, i.e., $\myAA(b,c)$ is of finite or affine type.
Then there exists a unique $\ZZ$-basis $\BB$ of $\myAA(b,c)$ such
that the semiring of positive elements of $\myAA(b,c)$ consists
precisely of positive integer linear combinations of elements of~$\BB$.
\end{theorem}

The uniqueness of $\BB$ is clear: it consists of all
indecomposable positive elements (those that cannot be written as
a sum of two positive elements).
We call $\BB$ the \emph{canonical} basis in $\myAA(b,c)$.

Our next result provides a parametrization of~$\BB$.
Let $Q = \ZZ^2$ be a lattice of rank 2 with a fixed basis
$\{\alpha_1,\alpha_2\}$.
We will sometimes write a point $\alpha = a_1 \alpha_1 +
a_2 \alpha_2 \in Q$ simply as $\alpha = (a_1, a_2)$.

\begin{theorem}
\label{th:canonicalbasis-parametrization}
In the situation of
Theorem~{\rm \ref{th:canonicalbasis}},
for every $\alpha = (a_1, a_2) \in Q$,
there is a unique basis element $x[\alpha] \in \BB$
of the form
\begin{equation}
\label{eq:leading-terms}
x[\alpha] = \frac{N_\alpha(y_1,y_2)}{y_1^{a_1}y_2^{a_2}},
\end{equation}
where $N_\alpha$ is a polynomial with constant term~$1$.
The correspondence
$\alpha \mapsto x[\alpha]$ is a bijection between $Q$ and $\BB$.
\end{theorem}

By its definition in Theorem~\ref{th:canonicalbasis}, the canonical basis $\BB$ is invariant under
any automorphism of $\myAA (b,c)$ that preserves the semiring of
positive elements; we refer to such an automorphism as
\emph{positive}.
In view of the exchange relations \eqref{eq:clusterrelations}, for
every $p \in \ZZ$, there is a positive automorphism
$\sigma_p$ of $\myAA(b,c)$ (for arbitrary $b$ and $c$)
acting on cluster variables by an involutive permutation $\sigma_p (y_m) = y_{2p-m}$.
The group of automorphisms of $\myAA(b,c)$ generated by all
$\sigma_p$ acts transitively on the set of clusters; it is easy to
see that this group is generated by any two elements $\sigma_p$ and
$\sigma_{p+1}$, for instance, by $\sigma_1$ and $\sigma_2$.

\begin{theorem}
\label{th:canonicalbasis-symmetry}
The bijection between $\BB$ and $Q$ given in
Theorem~{\rm \ref{th:canonicalbasis-parametrization}}
translates the action of each $\sigma_p$ on $\BB$ into a piecewise
linear transformation of~$Q$.
In particular, the action of $\sigma_1$ and $\sigma_2$ on $Q$ is given by
\begin{equation}
\label{eq:sigma-Q}
\sigma_1 (a_1, a_2) = (a_1, c \max (a_1, 0) - a_2), \quad
\sigma_2 (a_1, a_2) = (b \max (a_2, 0) - a_1, a_2) \ .
\end{equation}
\end{theorem}

\begin{remark}
\label{rem:taupm}
The piecewise linear transformations in \eqref{eq:sigma-Q} are a
special case of the transformations $\tau_+$ and $\tau_-$
appearing in \cite[(1.7)]{yga}.
Note that they can be viewed as ``tropicalizations" of
exchange relations \eqref{eq:clusterrelations}.
We expect that a similar interpretation exists for cluster
algebras of higher rank.
\end{remark}

\subsection{Explicit realization of $\BB$}
We now describe the canonical basis $\BB$ explicitly.
An element of the form $y_m^p y_{m+1}^q$ for some
$m\in\ZZ$ and $p,q \in \ZZ_{\geq 0}$ is called a \emph{cluster monomial}.

\begin{theorem}
\label{th:canonicalbasis-finite}
If $bc \leq 3$, i.e., $\myAA(b,c)$ is of finite type,
then the basis $\BB$ in Theorem~{\rm \ref{th:canonicalbasis}}
is the set of all cluster monomials.
\end{theorem}

Turning to the affine types, we notice that there are only two
cases to consider: $(b,c) = (2,2)$, or $(b,c) = (1,4)$
(the case $(b,c) = (4,1)$ reduces to the latter one by shifting
the indices of all cluster variables by~$1$).
In each of these two cases, we introduce an element $z \in
\myAA(b,c)$ by setting
\begin{equation}
\label{eq:z}
z =
\begin{cases}
y_0 y_3 - y_1 y_2 & \text{if $(b,c) = (2,2)$;} \\[.05in]
y_0^2 y_3 - (y_1+2) y_2^2 & \text{if $(b,c) = (1,4)$.}
\end{cases}
\end{equation}
Let $T_0, T_1, \dots$ be the sequence of Chebyshev polynomials of
the first kind given by $T_0 = 1$, and
$T_n (t + t^{-1}) = t^n + t^{-n}$ for $n > 0$.
In each of the affine types, we define the sequence $z_1, z_2,
\dots$ of elements of $\myAA(b,c)$ by setting $z_n = T_n (z)$.

\begin{theorem}
\label{th:canonicalbasis-affine}
If $(b,c)$ is one of the pairs $(2,2)$ or $(1,4)$
then the basis $\BB$ in Theorem~{\rm \ref{th:canonicalbasis}}
is the disjoint union of the set of all cluster monomials
and the set $\{z_n: n \geq 1\}$.
\end{theorem}

\begin{remark}
Theorems~\ref{th:canonicalbasis-finite} and \ref{th:canonicalbasis-affine}
imply in particular that in the finite or affine type,
every cluster variable $y_m$ is a positive element of $\myAA(b,c)$
in the sense of Definition~\ref{def:positivity}.
Jim Propp informed us that he (resp.~Gregg Musiker) found a combinatorial
proof of this fact for
$\myAA(2,2)$ (resp.~$\myAA(1,4)$).
Dylan Thurston told us that he found a topological
interpretation of $\myAA(2,2)$ leading to yet another proof of positivity
for this case.
As conjectured in \cite{fz-ClusterI}, the positivity of cluster variables is expected to hold
for cluster algebras of arbitrary rank; cf.~also Remark~\ref{rem:YZ-positivity} below.
\end{remark}

Following \cite{fz-ClusterI},
for $b, c$ arbitrary, we identify $Q$ with the root lattice associated with the Cartan
matrix (\ref{eq:cartan-rank2}), so that $\alpha_1$ and $\alpha_2$
become \emph{simple roots}
(see Section~\ref{sec:general-rank2}
for the background on the corresponding root system).
For the initial cluster variables $y_1$ and $y_2$, the correspondence
(\ref{eq:leading-terms}) takes the form
$$y_1 = \frac{1}{y_1^{-1}} = x[- \alpha_1],
\quad y_2 = \frac{1}{y_2^{-1}} = x[- \alpha_2].$$
As shown in \cite[Theorem~6.1]{fz-ClusterI}, the rest of the
cluster variables have the following description.

\begin{proposition}
\label{pr:ym-real-roots}
In a cluster algebra $\myAA(b,c)$ with $b,c$ arbitrary
(not necessarily of finite or affine type), every cluster variable
$y_m$ different from $y_1$ and $y_2$ has the form {\rm
(\ref{eq:leading-terms})}
for a positive real root $\alpha$.
This correspondence is a bijection between
the set of all positive real
roots associated with the Cartan matrix $A(b,c)$,
and the set of cluster variables
different from $y_1$ and $y_2$ in $\myAA(b,c)$.
\end{proposition}

In the affine type case $bc = 4$,
positive imaginary roots are all positive integer multiples of
the root $\delta$ given by
\begin{equation}
\label{eq:delta}
\delta =
\begin{cases}
\alpha_1 + \alpha_2 & \text{if $(b,c) = (2,2)$;} \\[.05in]
\alpha_1 + 2\alpha_2 & \text{if $(b,c) = (1,4)$.}
\end{cases}
\end{equation}
The following result complements
Proposition~\ref{pr:ym-real-roots}.

\begin{proposition}
\label{pr:z-imaginary}
Under the correspondence in
Theorem~\ref{th:canonicalbasis-parametrization}, we have
$x[n \delta] = z_n$ for all $n > 0$.
\end{proposition}

\begin{remark}
\label{rem:z-invariance}
As a consequence of
Theorem~\ref{th:canonicalbasis-symmetry} and
Proposition~\ref{pr:z-imaginary}, each $z_n$ is
invariant under all automorphisms $\sigma_p$ (see Proposition
\ref{pr:z-invariance} ).
In particular, this means that in the
``non-canonical" definition \eqref{eq:z}, the four
cluster variables $y_0, y_1, y_2$, and $y_3$ can be replaced by
any four consecutive cluster variables (see \eqref{eq:explicit-z}).
\end{remark}

\begin{remark}
\label{rem:leclerc}
In view of Theorem~\ref{th:canonicalbasis-affine}, no non-trivial
monomial in the $z_n$ belongs to the canonical basis $\BB$.
This is in contrast with the behavior of cluster variables, since
$\BB$ contains all cluster monomials.
This contrast between the ``real" generators $y_m$ and the
``imaginary" generators $z_n$ seems to be a special case of a
general phenomenon taking place in any cluster algebra of
infinite type (cf.~\cite{leclerc}).
\end{remark}

The rest of the paper is organized as follows.
In Section~\ref{sec:general-rank2}, after a necessary background on root systems of rank $2$,
we state and prove some results about cluster monomials which are valid in an arbitrary
algebra $\myAA (b,c)$ (not assumed to be of finite or affine type).
In particular, we show (Corollary~\ref{cor:linear-independence})
that the cluster monomials are always linearly independent.
Sharpening Proposition~\ref{pr:ym-real-roots}, we compute the Newton
polygon of every cluster variable $y_m$ viewed  as a Laurent
polynomial in $y_1$ and $y_2$ (Proposition~\ref{pr:triangleprop});
this description turns out to be crucial for the proof of
Theorem~\ref{th:canonicalbasis}.

Section~\ref{sec:proofs-finite} contains the proofs of our main
results (Theorems~\ref{th:canonicalbasis},
\ref{th:canonicalbasis-parametrization},
\ref{th:canonicalbasis-symmetry}, and
\ref{th:canonicalbasis-finite}) in the finite type case, while
Section~\ref{sec:affine-proofs} treats the affine types.
To avoid a case-by-case analysis of the algebras $\myAA(2,2)$ and $\myAA(1,4)$,
we introduce in Section~\ref{sec:affine-proofs} an interesting
$2$-parameter deformation of $\myAA(2,2)$, and extend to it most
of the results in question.

Finally, in Section~\ref{sec:coeffs-removed} we deal with more
general cluster algebras of rank~$2$, in which the two monomials
on the right hand side of \eqref{eq:clusterrelations} have
non-trivial coefficients.
We show (Proposition~\ref{pr:coeffs-removed} and Theorem~\ref{th:canonicalbasis-coeff})
that this seemingly more general case reduces to the coefficient-free case by
a rescaling of cluster variables.

\section{Cluster monomials and their Newton polygons}
\label{sec:general-rank2}

\subsection{Background on rank~$2$ root systems.}
In this section we work with an arbitrary cluster algebra
of rank~$2$; so the indefinite case $bc > 4$ is not excluded.
We start by recalling some basic facts about the root system associated with
a $2 \times 2$ Cartan matrix $A = A(b,c)$ (more details can be
found in~\cite{kac}).

First of all, the \emph{Weyl group} $W = W(A)$ is a group of linear
transformations of
$Q$ generated by two \emph{simple reflections} $s_1$ and $s_2$
whose action in the basis of simple roots is given by
\begin{equation}
\label{eq:s1-s2}
s_1 = \mat{-1}{b}{0}{1}\ ,\quad
s_2 = \mat{1}{0}{c}{-1}\ .
\end{equation}
Since both $s_1$ and $s_2$ are involutions, each element of $W$
is one of the following:
$$w_1 (m) = s_1 s_2 s_1 \cdots s_{\rem{m}}, \,\,
w_2 (m) = s_2 s_1 s_2 \cdots s_{\rem{m+1}} \ ;$$
here $\rem{m}$ stands for the element
of $\{1,2\}$ congruent to $m$ modulo $2$,
and both products are of length $m \geq 0$.
It is well known that $W$ is finite if and only if $bc \leq 3$.
The \emph{Coxeter number} $h$ of $W$ is the order of $s_1 s_2$ in
$W$; it is given by Table~\ref{tab:Coxeter number}.
In the finite case, $W$ is the dihedral group of order $2h$, and
its elements can be listed as follows: $w_1 (0) = w_2 (0) = e$
(the identity element), $w_1 (h) = w_2 (h) = w_0$ (the longest
element), and $2h - 2$ distinct elements $w_1 (m), w_2 (m)$
for $0 < m < h$.
In the infinite case, all elements $w_1 (m)$ and $w_2 (m)$ for $m > 0$
are distinct.

A vector $\alpha \in Q$ is a \emph{real root} for $A$ if
it is $W$-conjugate to a simple root.
Let $\Phi^{\rm re}$ denote the set of real roots for $A$.
It is known that
$\Phi^{\rm re} = \Phi_+^{\rm re} \cup (- \Phi_+^{\rm re})$,
where
$$\Phi_+^{\rm re} =
\{\alpha = a_1 \alpha_1 + a_2 \alpha_2 \in \Phi^{\rm re} : \ a_1, a_2 \geq
0\}$$
is the set of positive real roots.

In the finite case, $\Phi_+^{\rm re}$ has cardinality $h$,
and we have
$$\Phi_+^{\rm re} = \{w_1 (m) \alpha_{\rem{m+1}}: 0 \leq m < h\}
%= \{w_2 (m) \alpha_{\rem{m+2}}: 0 \leq m < h\}
\ .$$
In the infinite case, we have
$$\Phi_+^{\rm re} = \{w_i (m) \alpha_{\rem{m+i}}
\ : \ i \in \{1,2\}, \, m \geq 0\} \ ,$$
with all the roots $w_i (m) \alpha_{\rem{m+i}}$ distinct.

To introduce imaginary roots, consider the $W$-invariant
symmetric scalar product on~$Q$ given by
\begin{equation}
\label{eq:scalar-product}
(\alpha, \alpha) =
c a_1^2 - bc a_1 a_2 + b a_2^2 \quad
(\alpha = a_1 \alpha_1 + a_2 \alpha_2 \in Q).
\end{equation}
According to \cite{kac}, the set of imaginary
roots can be defined as follows:
\begin{equation}
\label{eq:imaginary-roots}
\Phi^{\rm im} = \{\alpha \in Q - \{0\} \, : \,
(\alpha, \alpha) \leq 0\}.
\end{equation}
Similarly to the case of real roots, we have
$\Phi^{\rm im} = \Phi_+^{\rm im} \cup (- \Phi_+^{\rm im})$,
where
$$\Phi_+^{\rm im} =
\{\alpha = a_1 \alpha_1 + a_2 \alpha_2 \in \Phi^{\rm im} : \ a_1, a_2 \geq
0\}$$
is the set of positive imaginary roots.
In view of (\ref{eq:imaginary-roots}) and
(\ref{eq:scalar-product}), a vector $a_1 \alpha_1 + a_2 \alpha_2$
from $Q$ belongs to $\Phi_+^{\rm im}$
if and only if $a_1, a_2 > 0$, and
\begin{equation}
\label{eq:imaginary-inequalities}
\frac{bc - \sqrt{bc(bc-4)}}{2b} \leq \frac{a_2}{a_1} \leq
\frac{bc + \sqrt{bc(bc-4)}}{2b} \ .
\end{equation}
In particular, in the affine type case $bc = 4$, positive imaginary roots
are precisely the
positive integer multiples of the minimal such root $\delta$ given by
(\ref{eq:delta}).
Note also that in the indefinite case $bc > 4$, the upper and
lower bounds for $a_2/a_1$ in (\ref{eq:imaginary-inequalities})
are irrational numbers (indeed, in this case,
$(bc - 3)^2 < bc(bc-4) < (bc-2)^2$, hence $bc(bc-4)$ is not a
perfect square); so both inequalities in \eqref{eq:imaginary-inequalities}
can be replaced by strict ones.

\subsection{Parameterizing cluster monomials}
Returning to cluster algebras, we start with the following result.

\begin{proposition}
\label{pr:cluster-monomials-leading-terms}
Every cluster monomial has a (unique) presentation of the form
{\rm (\ref{eq:leading-terms})}, with
$\alpha \in Q - \Phi^{\rm im}_+$.
This correspondence is a bijection between the set of all cluster monomials
and $Q - \Phi^{\rm im}_+$.
\end{proposition}

\begin{proof}
According to~\cite[Theorem~6.1]{fz-ClusterI} or
Proposition~\ref{pr:ym-real-roots} above,
every cluster variable has the form %(\ref{eq:leading-terms}), i.e.,
$y_m = x[\alpha(m)]$; here $\alpha(m) = - \alpha_m$ for
$m \in \{1,2\}$, and the rest of the cluster variables are in a
bijection with the set of all positive real roots.
Clearly, every cluster monomial $y_m^p y_{m+1}^q$ also has the
form (\ref{eq:leading-terms}), i.e., can be written as $x[\alpha]$
with $\alpha = p \alpha(m) + q \alpha(m+1)$.
To complete the proof it is enough to show the following:
\begin{eqnarray}
\label{eq:clusters-unimodular}
&&\text{For every $m \in \ZZ$, the vectors $\alpha(m)$ and $\alpha(m+1)$
form a
$\ZZ$-basis of $Q$.}\\
\label{eq:not-overlap}
&&\text{For every $m \in \ZZ$, the vectors $\alpha(m)$ and
$\alpha(m+1)$ are the only positive}\\
\nonumber
&&\text{real roots contained in the
additive semigroup $\ZZ_{\geq 0} \alpha(m) + \ZZ_{\geq 0} \alpha(m+1)$.}\\
%\hspace{-.2in}
\label{eq:cluster-union}
&&\text{The union
$\bigcup_{m \in \ZZ} (\ZZ_{\geq 0} \alpha(m) + \ZZ_{\geq 0} \alpha(m+1))$
is equal to $Q - \Phi^{\rm im}_+$.}
\end{eqnarray}

In the finite type, the properties
\eqref{eq:clusters-unimodular}--\eqref{eq:cluster-union}
are seen by an easy inspection (they are in fact special cases of the
results
in~\cite[Theorems~1.8,1.10]{yga} valid for all cluster algebras of finite
type).
So let us assume we are in the infinite type case $bc \geq 4$.
As shown in~\cite[Theorem~6.1]{fz-ClusterI}, the positive roots
$\alpha(m)$ are given by
\begin{equation}
\label{eq:denominator-roots-infinite}
\alpha (m+3) =  w_1 (m) \alpha_{\rem{m+1}}, \quad
\alpha (-m) =  w_2 (m) \alpha_{\rem{m+2}} \quad \quad
(m \geq 0) \ .
\end{equation}
In particular, we have $\alpha(3) = \alpha_1$ and
$\alpha(0) = \alpha_2$.
Thus, properties \eqref{eq:clusters-unimodular} and \eqref{eq:not-overlap}
hold for $m \in \{0,1,2\}$, and we also have
$$\bigcup_{m =0}^2 (\ZZ_{\geq 0} \alpha(m) + \ZZ_{\geq 0} \alpha(m+1))
= Q - (\ZZ_{> 0} \alpha_1 + \ZZ_{> 0} \alpha_2).$$

We now check  \eqref{eq:clusters-unimodular} and \eqref{eq:not-overlap}
for $m \geq 3$ (the case $m \leq 0$ follows by
the obvious symmetry between $\alpha_1$ and $\alpha_2$).
As shown in~\cite[(6.7)]{fz-ClusterI}, the roots $\alpha(m)$ for
$m \in \ZZ - \{1,2\}$ satisfy the recurrence
\begin{equation}
\label{eq:alpha-recurrence}
\alpha(m+1) + \alpha(m-1) =
\begin{cases}
b \alpha (m) & \text{if $m$ is odd;} \\[.15in]
c \alpha (m) & \text{if $m$ is even.}
\end{cases}
\end{equation}
Writing each $\alpha(m)$ as
$a_{m1} \alpha_1 + a_{m2} \alpha_2$, we conclude from
(\ref{eq:alpha-recurrence}) that
\begin{equation}
\label{eq:det-1}
\det \mat{a_{m1}}{a_{m+1,1}}{a_{m2}}{a_{m+1,2}}
= \det \mat{a_{m-1,1}}{a_{m1}}{a_{m-1,2}}{a_{m2}}
= \cdots = \det \mat{a_{21}}{a_{31}}{a_{22}}{a_{32}} = 1
\end{equation}
for $m \geq 2$.
This proves \eqref{eq:clusters-unimodular}, as well as the fact that for $m
\geq 3$, the
vectors $\alpha(m-1)$ and $\alpha(m+1)$ lie on the opposite sides
of the ray through $\alpha(m)$.
It follows that the sequence $(a_{m2}/a_{m1})_{m \geq 3}$ is
strictly increasing.
In view of (\ref{eq:imaginary-inequalities}), to finish the
proof of \eqref{eq:not-overlap} and \eqref{eq:cluster-union},
it is enough to show that
\begin{equation}
\label{eq:limit}
\lim_{m \to \infty}  \frac{a_{m2}}{a_{m1}} =
\frac{bc - \sqrt{bc(bc-4)}}{2b} \ .
\end{equation}
Let us abbreviate $u_m = a_{m2}/a_{m1}$.
%Since the sequence $(u_m)_{m \geq 3}$ is
%increasing, the limit in (\ref{eq:limit}) exists.
By (\ref{eq:denominator-roots-infinite}), for $m \geq 3$ we have
$\alpha(m+2) = s_1 s_2 \alpha(m)$.
Remembering (\ref{eq:s1-s2}), we conclude that
$$u_{m+2} = \frac{c-u_m}{bc-1 -bu_m} \quad (m \geq 3).$$
An elementary check shows that the transformation
$$u \mapsto f(u) =  \frac{c-u}{bc-1 -bu}
%= \frac{1}{b} \left(1 +  \frac{1}{bc-1 -bu}\right)
$$
preserves the interval $[0, \frac{bc - \sqrt{bc(bc-4)}}{2b}]$,
and we have $f(u) > u$ for $0 \leq u < \frac{bc - \sqrt{bc(bc-4)}}{2b}$.
This implies (\ref{eq:limit}) and completes the proof of
Proposition~\ref{pr:cluster-monomials-leading-terms}.
\end{proof}

\begin{figure}[ht]
\psset{unit=40pt}
\psset{labelsep=4pt}
\newgray{lightgray}{0.75}
\begin{pspicture}(-2.5,-2.5)(4.5,4.5)
\pspolygon[linestyle=none,fillcolor=lightgray,fillstyle=solid](0,0)(4.5,0)(4.5,2.25)

\psline[arrows=<->](-2.5,0)(4.5,0)
\psline[arrows=<->](0,-2.5)(0,4.5)
\multirput(-2,-2)(1,0){7}{\multirput(0,0)(0,1){7}{\psdots[dotscale=1](0,0)}}
\psdots[dotscale=1.5](-1,0)(0,-1)(1,0)(2,1)(3,2)(4,3)(0,1)(1,2)(2,3)(3,4)
\psdots[dotscale=1.5](1,1)(2,2)(3,3)(4,4)

\uput[135](1,1){$z_1$}
\uput[135](2,2){$z_2$}
\uput[135](3,3){$z_3$}
\uput[135](4,4){$z_4$}
\uput[135](-1,0){$y_1$}
\uput[135](0,-1){$y_2$}
\uput[135](1,0){$y_3$}
\uput[135](2,1){$y_4$}
\uput[135](3,2){$y_5$}
\uput[135](4,3){$y_6$}
\uput[135](0,1){$y_{0}$}
\uput[135](1,2){$y_{-1}$}
\uput[135](2,3){$y_{-2}$}
\uput[135](3,4){$y_{-3}$}
\uput[270](4.4,0){$\alpha_1$}
\uput[0](0,4.4){$\alpha_2$}
\end{pspicture}
\caption{Root lattice for $b=c=2$.  Lattice points in the closed shaded
  sector correspond to cluster monomials $y_3^py_4^q$.}
\label{fig:22rootlattice}
\end{figure}
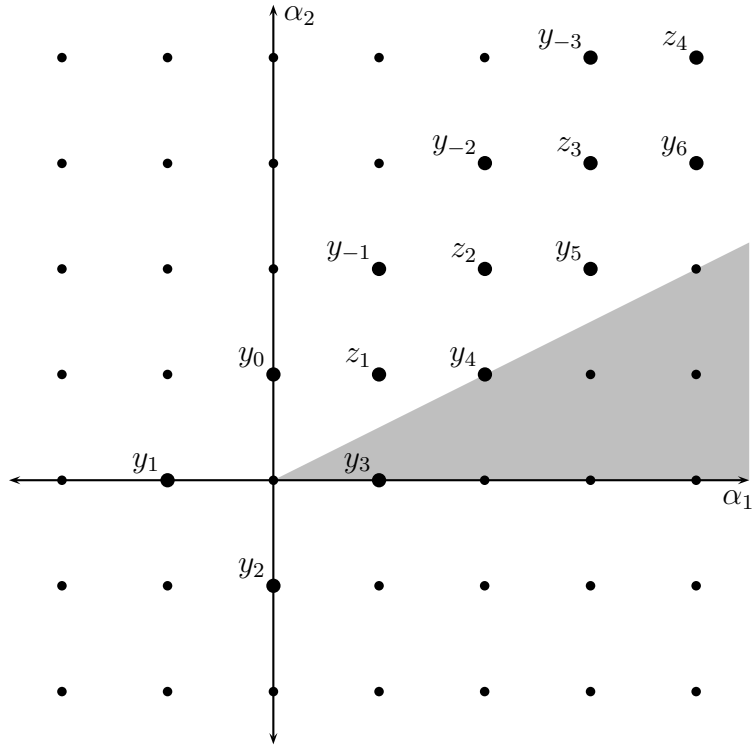

%\begin{figure}[ht]
%\begin{center}
%\setlength{\unitlength}{40pt}
%\begin{picture}(8,8)(-3,-3)
%\multiput(-2,-2)(0,1){7}{\multiput(0,0)(1,0){7}{\circle*{0.1}}}
%\put(-2.5,0){\line(1,0){7}}
%\put(0,-2.5){\line(0,1){7}}

%\put(-1,0){\circle*{0.2}}
%\put(-1,0){\makebox(0,0){$y_1$}}
%\put(0,-1){\circle*{0.2}}
%\put(0,-1){\makebox(0,0){$y_2$}}
%\put(1,0){\circle*{0.2}}
%\put(1,0){\makebox(0,0){$y_3$}}
%\put(2,1){\circle*{0.2}}
%\put(2,1){\makebox(0,0){$y_4$}}
%\put(3,2){\circle*{0.2}}
%\put(3,2){\makebox(0,0){$y_5$}}
%\put(4,3){\circle*{0.2}}
%\put(4,3){\makebox(0,0){$y_6$}}
%\put(0,1){\circle*{0.2}}
%\put(0,1){\makebox(0,0){$y_0$}}
%\put(1,2){\circle*{0.2}}
%\put(1,2){\makebox(0,0){$y_{-1}$}}
%\put(2,3){\circle*{0.2}}
%\put(2,3){\makebox(0,0){$y_{-2}$}}
%\put(3,4){\circle*{0.2}}
%\put(3,4){\makebox(0,0){$y_{-3}$}}

%\put(1,1){\circle*{0.2}}
%\put(1,1){\makebox(0,0){$z_1$}}
%\put(2,2){\circle*{0.2}}
%\put(2,2){\makebox(0,0){$z_2$}}
%\put(3,3){\circle*{0.2}}
%\put(3,3){\makebox(0,0){$z_3$}}
%\put(4,4){\circle*{0.2}}
%\put(4,4){\makebox(0,0){$z_4$}}

%\put(2,0){\line(-4,3){0.8}}

%  \qbezier(32,19)(24,25)(16,31)

%\end{picture}
%\end{center}
%\caption{Root lattice for the affine type $b=c=2$}
%\label{fig:22rootlattice}
%\end{figure}

\begin{example}
We illustrate Proposition~\ref{pr:cluster-monomials-leading-terms}
by explicitly computing the family $(\alpha(m))_{m \in \ZZ}$ for the
affine types. First, let $(b,c) = (2,2)$.
Using the recurrence (\ref{eq:alpha-recurrence}) and
induction on~$|m|$, we obtain that
\begin{equation}
\label{eq:alpham-22}
\alpha(m)  =
\begin{cases}
(m-2) \alpha_1 + (m-3) \alpha_2 & \text{if $m \geq 2$;} \\[.15in]
-m \alpha_1 + (1-m) \alpha_2 & \text{if $m \leq 1$.}
\end{cases}
\end{equation}
These roots and the corresponding cluster variables are shown in
Figure~\ref{fig:22rootlattice}.

Similarly, for $(b,c) = (1,4)$, we get
\begin{equation}
\label{eq:alpham-14}
\rem{m} \alpha(m)  =
\begin{cases}
(m-2) \alpha_1 + 2(m-3) \alpha_2 & \text{if $m \geq 2$;} \\[.15in]
-m \alpha_1 + 2 (1-m) \alpha_2 & \text{if $m \leq 1$.}
\end{cases}
\end{equation}
\end{example}

\begin{corollary}
\label{cor:linear-independence}
In any cluster algebra of rank~$2$, the cluster monomials are
linearly independent.
\end{corollary}

\begin{proof}
For $\gamma = g_1 \alpha_1 + g_2 \alpha_2 \in Q$, we abbreviate
$y^\gamma = y_1^{g_1} y_2^{g_2}$.
We will use the product partial order on $Q = \ZZ^2$:
$$\gamma_1 \geq \gamma_2 \,\, \Leftrightarrow \,\,
\gamma_1 - \gamma_2 \in Q_+ = \ZZ_{\geq 0} \alpha_1
+ \ZZ_{\geq 0} \alpha_2 \ .$$
By Proposition~\ref{pr:cluster-monomials-leading-terms},
cluster monomials can be parametrized by $Q - \Phi^{\rm im}_+$
so that the cluster monomial corresponding to
$\alpha \in Q - \Phi^{\rm im}_+$ has the form
\begin{equation}
\label{eq:cluster-monomial-Laurent}
x[\alpha] = y^{- \alpha} + \sum_{\gamma > -\alpha} c_\gamma
y^\gamma \ .
\end{equation}
Now suppose that a (finite) integer linear combination
of the cluster monomials $x[\alpha]$ is equal to~$0$.
Let $S \subset Q - \Phi^{\rm im}_+$ be the set of all $\alpha$
such that $x[\alpha]$ occurs with a non-zero coefficient in this
linear combination.
If $S$ is non-empty, pick a maximal element $\alpha \in S$;
in view of (\ref{eq:cluster-monomial-Laurent}), the (Laurent)
monomial $y^{- \alpha}$ does not occur in any $x[\beta]$ for
$\beta \in S - \{\alpha\}$, which gives a desired contradiction.
\end{proof}

We now turn our attention to positive automorphisms $\sigma_p$ of
$\myAA (b,c)$ which appear in Theorem~\ref{th:canonicalbasis-symmetry}.
By definition, each $\sigma_p$ acts on the set of all cluster
monomials; identifying the latter set with $Q - \Phi^{\rm im}_+$
as in Proposition~\ref{pr:cluster-monomials-leading-terms}, we
obtain the action of $\sigma_p$ on $Q - \Phi^{\rm im}_+$.

\begin{proposition}
\label{pr:sigma-cluster-monomials}
The action of $\sigma_1$ and $\sigma_2$ on $Q - \Phi^{\rm im}_+$
is given by {\rm \eqref{eq:sigma-Q}}.
\end{proposition}

\begin{proof}
By symmetry, it is enough to prove the claim for $\sigma_1$.
In the finite type case, it is seen by inspection, so we shall
deal with  the infinite type case $bc \geq 4$.
Remembering the definition, we see that the action of $\sigma_1$
on $Q - \Phi^{\rm im}_+$ is given as follows:
$\sigma_1 (\alpha(m)) = \alpha(2-m)$ for $m \in \ZZ$,
and $\sigma_1$ acts linearly in each cone
$\ZZ_{\geq 0} \alpha(m) + \ZZ_{\geq 0} \alpha(m+1)$; here
$\alpha (1) = - \alpha_1, \alpha (2) = - \alpha_2$, and the rest
of the $\alpha(m)$ are all positive real roots labeled according
to \eqref{eq:denominator-roots-infinite}.
This description implies that $\sigma_1 (-\alpha_1) = - \alpha_1$,
and $\sigma_1 (\alpha (m)) = s_2 \alpha (m)$ for $m \neq 1$, i.e.,
when $\alpha (m)$ is either a positive real root, or $- \alpha_2$.
It follows that the action of $\sigma_1$ on
$\alpha = (a_1, a_2) \in Q - \Phi^{\rm im}_+$ is given as follows:
$$\sigma_1(\alpha) =
\begin{cases}
(a_1, - a_2) & \text{if $a_1 \leq 0$;} \\[.05in]
s_2 (a_1, a_2) = (a_1, ca_1- a_2) & \text{if $a_1 \geq 0$.}
\end{cases}
$$
This is clearly equivalent to the first equality in
\eqref{eq:sigma-Q}, and we are done.
\end{proof}

\subsection{Newton polygons of cluster variables}
Sharpening (\ref{eq:cluster-monomial-Laurent}), we now give
an explicit description of the Newton polygon of any cluster variable $y_m$.
Recall that the Newton polygon ${\rm Newt} (x)$ of a Laurent
polynomial $x \in \ZZ[y_1^{\pm 1}, y_2^{\pm 1}]$
is the convex hull in $Q_\RR = \RR \alpha_1 \oplus \RR \alpha_2$
 of all lattice points
$\gamma$ such that the monomial $y^\gamma$
appears with a non-zero coefficient in the Laurent expansion of~$x$.
We will say that $x$ is \emph{monic} if
every monomial corresponding to a vertex of ${\rm Newt}(x)$ appears
in the Laurent expansion of~$x$ with coefficient~$1$.
For $\alpha=a_1 \alpha_1 + a_2 \alpha_2 \in Q_+$, let
$\Delta(\alpha)$ denote the triangle (possibly degenerate) in $Q_\RR$
with vertices $- \alpha, - \alpha + b a_2 \alpha_1$, and
$- \alpha + c a_1 \alpha_2$.

\begin{proposition}
\label{pr:triangleprop}
For every positive real root $\alpha$,
the corresponding cluster variable $x[\alpha]$ is a monic
Laurent polynomial in $y_1$ and $y_2$,
and its Newton polygon ${\rm Newt}(x[\alpha])$ is
equal to the triangle $\Delta(\alpha)$.
\end{proposition}

\begin{proof}
%In the finite type, both claims follow by inspection.
%Thus, we assume that $bc \geq 4$, and so all the roots $\alpha(m)$
%for $m \in \ZZ$ are distinct.
Using the obvious symmetry, it is enough to prove the proposition
for the cluster variables $y_m$ with $m \geq 3$.
The cluster relations \ref{eq:clusterrelations} imply immediately that
$$y_3 = x[\alpha_1] = y^{-\alpha_1} + y^{-\alpha_1+ c \alpha_2}$$
and
$$y_4 = x[b \alpha_1 + \alpha_2] = y^{-b \alpha_1 - \alpha_2} (y^{c
\alpha_2}
+ 1)^b + y^{-\alpha_2} \ ;$$
it follows that both $y_3$ and $y_4$ satisfy the proposition.
Now let $m \geq 4$ be such that each of the cluster variables
$y_{m-1}, y_m$ and $y_{m+1}$ is different from $y_1$ and $y_2$
(this condition is vacuous for the infinite types).
Proceeding by induction on $m$, we assume that the proposition
holds for $y_{m-1}$ and $y_m$, and will prove it for $y_{m+1}$.
Let $m$ be even (the case of $m$ odd is entirely similar).
To find ${\rm Newt}(y_{m+1})$, we compute the Newton polygons of
both sides of the relation $y_{m-1} y_{m+1} = y_m^c + 1$.
Since the Newton polygon of the product of two Laurent polynomials
is the Minkowski sum of their Newton polygons, we obtain

\begin{equation}
\label{eq:sum-polygons}
{\rm Newt}(x[\alpha(m-1)]) + {\rm Newt}(x[\alpha(m+1)]) =
{\rm Newt}(x[\alpha(m)]^c + 1) \ .
\end{equation}

We will use the following well-known properties
of convex polygons in $\RR^2$ and their Minkowski
addition (see~\cite[Section~15.3]{grunbaum}).
With every convex $n$-gon $\Pi$ in $\RR^2$ there is associated a
set $V(\Pi)$  of $n$ non-zero vectors in $\RR^2$ defined as follows: each
side $AB$ of
$\Pi$ contributes a vector of the same length as $AB$ whose direction is
that of the
outward normal to $AB$.
In the degenerate case $n = 2$ when $\Pi$ is just a line segment,
the set $V(\Pi)$ is a pair of opposite vectors, each of length
equal to $\Pi$ and perpendicular to $\Pi$.
Clearly, the sum of $n$ vectors from $V(\Pi)$ is $0$, and no two of them
are positively proportional.
Conversely, every $n$-set of non-zero vectors in $\RR^2$ with
these properties has the form $V(\Pi)$, where a convex $n$-gon
$\Pi$ is determined uniquely up to translations.
Furthermore, for every two convex polygons $\Pi_1$ and $\Pi_2$ in
$\RR^2$, the set $V(\Pi_1 + \Pi_2)$ associated with their Minkowski sum
is obtained by taking the union of $V(\Pi_1)$ and $V(\Pi_2)$ (as a multiset)
and replacing any pair of positively proportional vectors by their sum.
This implies in particular the \emph{cancellation property} for the
Minkowski addition:
a convex polygon $\Pi_2$ is uniquely determined by $\Pi_1$ and
$\Pi_1 + \Pi_2$.

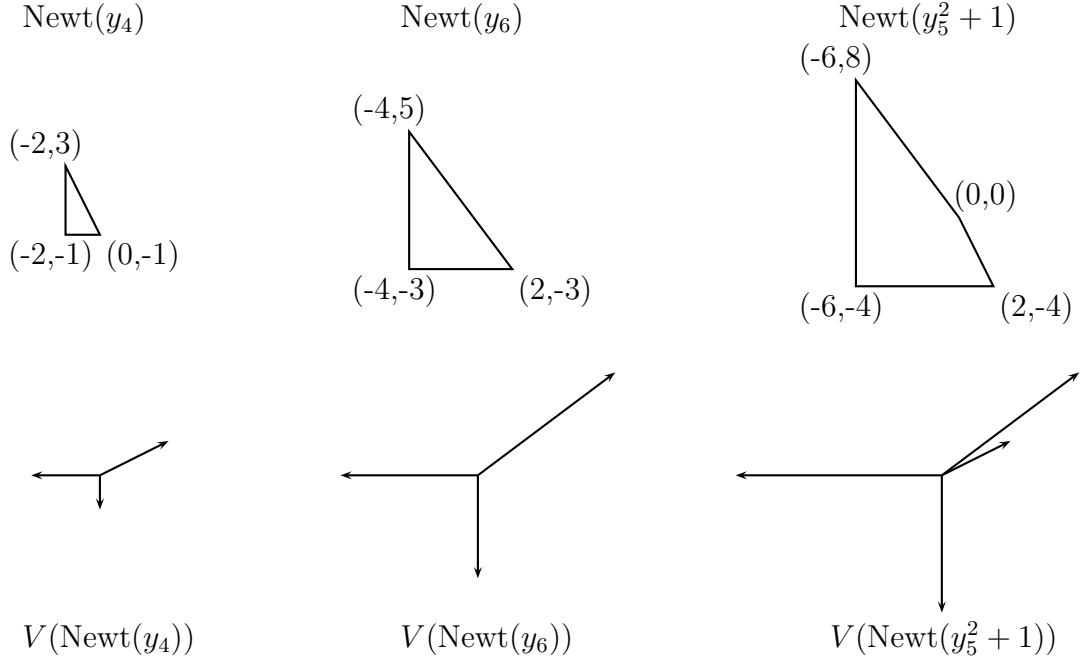
\begin{figure}[ht]
\begin{center}
\psset{unit=6.5pt}
\psset{labelsep=3pt}
\begin{pspicture}(0,0)(65,40)
\put(0,3){
 \put(7,24){
    \pspolygon(-2,-1)(0,-1)(-2,3)
    \uput[225](-2,-1){\put(-3,0){(-2,-1)}}
    \uput[315](0,-1){\put(0,0){(0,-1)}}
    \uput[135](-2,3){\put(-3,0){(-2,3)}}
    \put(-4.5,11.1){Newt$(y_4)$}
  }
  \put(7,9){
   \psline[arrows=->](0,0)(-4,0)
    \psline[arrows=->](0,0)(0,-2)
    \psline[arrows=->](0,0)(4,2)
    \put(-4.5,-10){$V$(Newt$(y_4)$)}
  }
}
\put(24,3){
  \put(5,24){
    \pspolygon(-4,-3)(2,-3)(-4,5)
    \uput[225](-4,-3){\put(-3,0){(-4,-3)}}
    \uput[315](2,-3){\put(0,0){(2,-3)}}
    \uput[135](-4,5){\put(-3,0){(-4,5)}}
    \put(-4.5,11.1){Newt$(y_6)$}
  }
  \put(5,9){
    \psline[arrows=->](0,0)(-8,0)
    \psline[arrows=->](0,0)(0,-6)
    \psline[arrows=->](0,0)(8,6)
    \put(-4.5,-10){$V$(Newt$(y_6)$)}
  }
}
\put(52,3){
  \put(5,24){
    \pspolygon(-6,-4)(2,-4)(0,0)(-6,8)
    \uput[225](-6,-4){\put(-3,0){(-6,-4)}}
    \uput[315](2,-4){\put(0,0){(2,-4)}}
    \uput[135](0,0){\put(0,0){(0,0)}}
    \uput[135](-6,8){\put(-3,0){(-6,8)}}
    \put(-7,11.1){Newt$(y_5^2+1)$}
  }
  \put(4,9){
    \psline[arrows=->](0,0)(4,2)
    \psline[arrows=->](0,0)(8,6)
   \psline[arrows=->](0,0)(0,-8)
    \psline[arrows=->](0,0)(-12,0)
    \put(-6.5,-10){$V$(Newt$(y_5^2+1)$)}
  }
}
\end{pspicture}
\caption{Minkowski addition for the exchange relation $y_4y_6=y_5^2+1$
in the case $(b,c)=(2,2)$}
\end{center}
\label{fig:minkowski}
\end{figure}

Applying the cancellation property to (\ref{eq:sum-polygons}) we see
that the desired equality
$${\rm Newt}(x[\alpha(m+1)]) = \Delta(\alpha(m+1))$$
becomes a consequence of
\begin{equation}
\label{eq:sum-triangles}
\Delta(\alpha(m-1)) + \Delta(\alpha(m+1)) =
{\rm Conv} (c \Delta(\alpha(m)) \cup \{0\}),
\end{equation}
where ${\rm Conv}$ stands for the convex hull in $Q_\RR = \RR^2$.

Our proof of (\ref{eq:sum-triangles}) is illustrated by Figure~\ref{fig:minkowski}.
We first compute the polygon on the right hand side.
Note that if $\alpha = a_1 \alpha_1 + a_2 \alpha_2 \in Q$ is such that
$a_1, a_2 > 0$ then the side of the triangle $\Delta(\alpha)$
opposite to the vertex $- \alpha$ lies on the straight line
$$\{g_1 \alpha_1 + g_2 \alpha_2 \, : \, ca_1g_1+ba_2g_2 +
(\alpha,\alpha) = 0\},$$
where $(\alpha,\alpha)$ is given by~(\ref{eq:scalar-product}).
If $\alpha$ is a positive real root then $(\alpha,\alpha) > 0$,
hence $\Delta(\alpha)$ does not contain the origin.
We conclude that ${\rm Conv} (c \Delta(\alpha(m)) \cup \{0\})$
is the quadrilateral with vertices $- c \alpha(m),
- c \alpha(m) + bc a_{m2} \alpha_1, 0$, and
$- c\alpha(m) + c^2 a_{m1} \alpha_2$, where as before, we write
$\alpha(m)= a_{m1} \alpha_1 + a_{m2} \alpha_2$.
It follows easily that the corresponding set
$V({\rm Conv} (c \Delta(\alpha(m)) \cup \{0\}))$
consists of the following four vectors written
in the basis ${\alpha_1,\alpha_2}$:
\begin{equation}
\label{eq:vset-1}
\{(-c^2 a_{m1},0), (0,-bc a_{m2}),
(c a_{m2}, bc a_{m2} -c a_{m1}),
(c^2 a_{m1} - c a_{m2}, c a_{m1})\} \ .
\end{equation}
The position of the quadrilateral itself is uniquely determined by this set
and the property that its minimal vertex (with respect to the
product order on $\RR^2$) is
$- c \alpha(m) = (-c a_{m1}, -c a_{m2})$.

Turning to the left hand side of \eqref{eq:sum-triangles}, we use
the following easily checked property: if
$\alpha=a_1 \alpha_1 + a_2 \alpha_2 \in Q_+$
then
$$V(\Delta(\alpha)) = \{(- c a_1,0), (0,- ba_2),
(c a_1, ba_2)\}$$
(with the convention that if, say, $a_1 > 0$ and $a_2 = 0$
then the zero vector $(0,- ba_2)$ does not appear in $V(\Delta(\alpha))$).
Using the equality $\alpha(m-1) + \alpha(m+1) = c \alpha (m)$
from \eqref{eq:alpha-recurrence}, we see that
$V(\Delta(\alpha(m-1)) + \Delta(\alpha(m+1)))$ is equal to
\begin{equation}
\label{eq:vset-2}
\{(-c^2 a_{m1},0), (0,-bc a_{m2}), (c a_{m-1,1}, ba_{m-1,2}),
(c a_{m+1,1}, ba_{m+1,2})\} \ .
\end{equation}
It follows that $\Delta(\alpha(m-1)) + \Delta(\alpha(m+1))$ is a
quadrilateral;
again, its position is uniquely determined by
$V(\Delta(\alpha(m-1)) + \Delta(\alpha(m+1)))$ and
the property that its minimal vertex  is
$- (\alpha(m-1) + \alpha(m+1)) = - c \alpha (m)$.

To show \eqref{eq:sum-triangles}, it remains to prove that the sets in
\eqref{eq:vset-1} and \eqref{eq:vset-2} are equal to each other.
Remembering that the vectors in each of the sets must sum to zero,
we see it suffices to show that
$$a_{m2} = a_{m-1,1}, \,\, c a_{m1} = ba_{m+1,2} \ .$$
Once again using \eqref{eq:alpha-recurrence}, we can rewrite these
equalities as
\begin{equation}
\label{eq:m-to-m+2}
a_{m+1,2} = c a_{m-1,1} - a_{m-1,2}, \quad
a_{m-1,1} = b a_{m+1,2} - a_{m+1,1},
\end{equation}
which is an immediate consequence of
$\alpha(m+1) = s_1 s_2 \alpha(m-1)$.
This concludes the proof of the equality
${\rm Newt}(x[\alpha(m+1)]) = \Delta(\alpha(m+1))$.

To complete the proof of Proposition~\ref{pr:triangleprop}, it
remains to show that if $y_{m-1}$ and $y_m$ are monic then so is $y_{m+1}$.
In view of the exchange relations, this is a consequence
of the following easily verified general property:
a Laurent polynomial equal to the ratio of
two monic Laurent polynomials is itself monic.
\end{proof}

We conclude this section with the following corollary of
Proposition~\ref{pr:triangleprop}.

\begin{proposition}
\label{pr:no-positive-terms}
Let $y = y_m^p y_{m+1}^q$ be a cluster monomial containing at
least one cluster variable different from $y_1$ and $y_2$.
Then there exists a non-zero linear form $\varphi_m$ on $Q_\RR$ such that
$$\varphi_m (g_1 \alpha_1 + g_2 \alpha_2) = c_{m1} g_1 + c_{m2}
g_2, \quad c_{m1}, c_{m2} \geq 0,$$
and ${\rm Newt}(y)$ lies in the open half-plane $\{\varphi_m < 0\}$.
In particular, ${\rm Newt}(y)$ has empty intersection with the positive
quadrant $Q_+ = \ZZ_{\geq 0} \alpha_1 + \ZZ_{\geq 0} \alpha_2$.
\end{proposition}

\begin{proof}
It is enough to consider the case $m \geq 2$;
by Proposition~\ref{pr:finite-type-periodicity},
in the finite type case we can also assume
without loss of generality that $m \leq \frac{h}{2} + 2$.
If $m = 2$ then $y = y_2^p y_3^q$ for some $p \geq 0$ and $q > 0$;
in this case, ${\rm Newt}(y)$ is a line segment parallel to
$\alpha_2$ and lying in the open half-plane
$\{g_1 \alpha_1 + g_2 \alpha_2 \ : \ g_1 < 0\}$.
Thus, the proposition holds for $m = 2$.

For $m \geq 3$, we define a linear form $\varphi_m$ on $Q_\RR$ by setting
$$\varphi_m (g_1 \alpha_1 + g_2 \alpha_2) =
(a_{m+1,2} + c a_{m1} - a_{m2}) g_1 +
(a_{m,1} + b a_{m+1,2} - a_{m+1,1}) g_2 \ .$$
Remembering \eqref{eq:m-to-m+2}, we can rewrite $\varphi_m$ as
\begin{equation}
\label{eq:phi}
\varphi_m (g_1 \alpha_1 + g_2 \alpha_2) =
(a_{m+2,2} + a_{m+1,2}) g_1 +
(a_{m-1,1} + a_{m,1}) g_2 \ ,
\end{equation}
so the coefficients of $g_1$ and $g_2$ in $\varphi_m$ are nonnegative.
Since
$${\rm Newt}(y) = p \Delta(\alpha(m)) + q \Delta(\alpha(m+1)),$$
to complete the proof it is enough to show that $\varphi_m$ takes negative
values
at all vertices of each of the triangles $\Delta(\alpha(m))$ and
$\Delta(\alpha(m+1))$.
Clearly, both $\varphi_m (- \alpha(m))$ and $\varphi_m (-\alpha(m+1))$
are negative since $\varphi_m$ has positive coefficients.
As for the remaining four vertices, a direct calculation using
\eqref{eq:det-1} yields
$$\varphi_m (- \alpha(m) + c a_{m1} \alpha_2) =
\varphi_m (- \alpha(m+1) + b a_{m+1,2} \alpha_1) = -1,$$
$$\varphi_m (- \alpha(m) + b a_{m2} \alpha_1) =
-(\alpha(m), \alpha(m)) - 1,$$
$$\varphi_m (- \alpha(m+1) + c a_{m+1,1} \alpha_2) =
-(\alpha(m+1), \alpha(m+1)) - 1.$$
Since $\alpha (m)$ and $\alpha(m+1)$ are real roots,
both $(\alpha(m), \alpha(m))$ and $(\alpha(m+1), \alpha(m+1))$ are positive.
Hence all the above values of $\varphi_m$ are negative, as desired.
\end{proof}

In the affine types, the linear forms $\varphi_m$ given by \eqref{eq:phi}
can be computed explicitly using \eqref{eq:alpham-22} and
\eqref{eq:alpham-14}.
This implies the following corollary which will be used in
Section~\ref{sec:affine-proofs}.

\begin{corollary}
\label{cor:phi-affine}
In an affine type, for $m \geq 3$, every Laurent monomial
$y_1^{g_1} y_2^{g_2}$ that occurs in the expansion of a non-trivial cluster
monomial in $y_m$ and $y_{m+1}$ satisfies the following condition:
\begin{itemize}

\item $(2m-3)g_1 + (2m-5)g_2 < 0$ if $(b,c) = (2,2)$.

\item $(3m-3 - \rem{m})g_1 + \frac{3m-6 -\rem{m}}{2}g_2 < 0$ if $(b,c) =
(1,4)$.
\end{itemize}
\end{corollary}

\section{Proofs for finite types}
\label{sec:proofs-finite}

Throughout this section, $\myAA(b,c)$ is a cluster algebra of finite type.
It is enough to consider three cases: type~$A_2$ with
$(b,c) = (1,1)$, type~$B_2$ with $(b,c) = (1,2)$, and type~$G_2$ with
$(b,c) = (1,3)$.
Note that in these cases there are no imaginary roots, so
Theorem~\ref{th:canonicalbasis-parametrization}
(resp.~Theorem~\ref{th:canonicalbasis-symmetry})  becomes a
consequence of Theorem~\ref{th:canonicalbasis-finite} and
Proposition~\ref{pr:cluster-monomials-leading-terms}
(resp.~Proposition~\ref{pr:sigma-cluster-monomials}).

To prove Theorems~\ref{th:canonicalbasis} and
\ref{th:canonicalbasis-finite}, it suffices to show the following:
\begin{eqnarray}
\label{eq:y-positivity-finite}
&&\text{Every cluster variable is  a positive element of
$\myAA(b,c)$.}\\
\label{eq:spanning-finite}
&&\text{Cluster monomials linearly span $\myAA(b,c)$.}\\
\label{eq:positive-converse-finite}
&&\text{If $y \in \myAA(b,c)$ is written as a linear combination of
cluster monomials then}
\\
\nonumber
&&\text{its every coefficient
is equal to some
coefficient in the Laurent polynomial}
\\
\nonumber
&&\text{expansion of $y$ with respect to some cluster $\{y_m,y_{m+1}\}$.}
\end{eqnarray}

Recall from Proposition~\ref{pr:finite-type-periodicity} that
every cluster variable $y_m$ is equal to one of the variables
$y_1, y_2, \dots, y_{h+2}$.
To prove \eqref{eq:y-positivity-finite}, is is enough to show that each of
the elements
$y_3, \dots, y_{h+2}$ is a positive Laurent polynomial in
$y_1$ and $y_2$ (since all the clusters are interchangeable).
This is clear from their explicit calculation (the formulas below
are obtained from those in the proof of \cite[Theorem~6.1]{fz-ClusterI}
by specializing all the coefficients to $1$):

\smallskip

\noindent  \textsl{Type $A_2$:} $(b,c) = (1,1)$.
\begin{equation}
\label{eq:ym-Laurent-A2}
y_3 = \frac{y_2+1}{y_1}, \quad y_4 =
\frac{y_1 + y_2 + 1}{y_1 y_2}, \quad y_5 = \frac{y_1 + 1}{y_2} \ .
\end{equation}

\smallskip

\noindent \textsl{Type $B_2$:} $(b,c) = (1,2)$.
\begin{equation}
\label{eq:ym-Laurent-B2}
y_3 = \frac{y_2^2 + 1}{y_1}, \quad
y_4 = \frac{y_1 + y_2^2 + 1}{y_1 y_2}, \quad
y_5 = \frac{(y_1 + 1)^2 + y_2^2}{y_1 y_2^2}, \quad
y_6 = \frac{y_1 + 1}{y_2}\ .
\end{equation}

\smallskip

\noindent  \textsl{Type $G_2$:} $(b,c) = (1,3)$.
\begin{eqnarray}
\label{eq:ym-Laurent-G2}
&& y_3 = \frac{y_2^3 + 1}{y_1}, \quad y_4 = \frac{y_1 + y_2^3 + 1}{y_1 y_2},
\quad y_5 = \frac{(y_1 + 1)^3 + y_2^3 (y_2^3 + 3y_1 + 2)}
{y_1^2 y_2^3},\\
\nonumber
&&y_6 = \frac{(y_1 + 1)^2 + y_2^3}{y_1 y_2^2}, \quad
y_7 = \frac{(y_1 + 1)^3 + y_2^3}{y_1 y_2^3},
\quad y_8 = \frac{y_1 + 1}{y_2}\ .
\end{eqnarray}

\smallskip

To prove \eqref{eq:spanning-finite}, we notice that cluster monomials are
the monomials
in the cluster variables $y_m$ that do not contain
products of the form $y_m y_{m+n}$ for $n \geq 2$.
Our first task is to produce the relations that express
every such ``forbidden" product as a
linear combination of cluster monomials.
Without loss of generality we can assume that $n \leq (h+2)/2$.
For $n = 2$, the relations in question are just the exchange
relations (\ref{eq:clusterrelations}).
It remains to treat $n = 3$ in the $B_2$ case,
and $n = 3,4$ in the $G_2$ case.
A direct calculation using \eqref{eq:ym-Laurent-B2} and
\eqref{eq:ym-Laurent-G2} yields the following:

\smallskip

\noindent
\textsl{Type~$B_2$:} $(b,c) = (1,2)$.
\begin{equation}
\label{eq:straightening-B2}
y_m y_{m+3} =
\begin{cases}
y_{m-1} +  y_{m+1} & \text{if $m$ is odd;} \\[.05in]
y_{m+2} +  y_{m+4} & \text{if $m$ is even.}
\end{cases}
\end{equation}

\noindent
\textsl{Type~$G_2$:} $(b,c) = (1,3)$.
\begin{equation}
\label{eq:straightening-G2-3}
y_m y_{m+3} =
\begin{cases}
y_{m-1} +  y_{m+1}^2 & \text{if $m$ is odd;} \\[.05in]
y_{m+2}^2 +  y_{m+4} & \text{if $m$ is even.}
\end{cases}
\end{equation}
\begin{equation}
\label{eq:straightening-G2-4}
y_m y_{m+4} =
\begin{cases}
y_{m-2} +  y_{m+2} + 3 & \text{if $m$ is odd;} \\[.05in]
y_{m-2} +  y_{m+2} & \text{if $m$ is even.}
\end{cases}
\end{equation}

\smallskip

To complete the proof of~\eqref{eq:spanning-finite}, it remains to show that
the above
relations form a system of \emph{straightening} relations, i.e.,
that their repeated application allows one to express every
monomial in the cluster variables as a linear combination of
cluster monomials.

In type $A_2$, the exchange relations are of the form
$y_m y_{m+2} = y_m + 1$.
Note that both terms on the right have smaller degree than the
product on the left hand side.
It follows that if we take any monomial $M$
in the cluster variables which has $y_m y_{m+2}$ as a factor, and
replace this factor by $y_m + 1$, then $M$ will be
written as a sum of two monomials of smaller degree.
The fact that $M$ is a linear combination of cluster
monomials follows at once by induction on the degree of $M$.

In types $B_2$ and $G_2$, the argument is similar.
The only difference is that we use the degree of
monomials in the cluster variables based on the following
grading: in type~$B_2$ (resp.~$G_2$) we set
${\rm deg} (y_m) = 3$ (resp.~$5$) for~$m$ odd,
and ${\rm deg} (y_m) = 2$ (resp.~$3$) for~$m$ even.
By inspection, every monomial on the right hand side of
\eqref{eq:straightening-B2}--\eqref{eq:straightening-G2-4}
(and of every exchange relation for $B_2$ and $G_2$) has smaller
degree than the ``forbidden" product on the left hand side of the same
relation.
Thus, the proof of \eqref{eq:spanning-finite} is completed by the same
argument.

To prove \eqref{eq:positive-converse-finite}, suppose that an element $y \in \myAA(b,c)$
is written as a linear combination of cluster monomials.
Again, since all the clusters are interchangeable, it is enough to show
that the coefficient of a cluster monomial $y_1^p y_2^q$ in this
linear combination is equal to the coefficient of $y_1^p y_2^q$ in the Laurent polynomial expansion
of $y$ with respect to $\{y_1,y_2\}$.
This follows at once from Proposition~\ref{pr:no-positive-terms},
which implies that, for $p, q \geq 0$, the monomial $y_1^p y_2^q$ does not occur in the
Laurent expansion of any other cluster monomial.
This concludes the proofs of Theorems~\ref{th:canonicalbasis} and
\ref{th:canonicalbasis-finite} for the finite type.

\section{Proofs for affine types}
\label{sec:affine-proofs}

In this section, we prove Theorems~\ref{th:canonicalbasis},
\ref{th:canonicalbasis-parametrization}, \ref{th:canonicalbasis-symmetry},
\ref{th:canonicalbasis-affine}, and
Proposition~\ref{pr:z-imaginary} for the algebras
$\myAA(2,2)$ and $\myAA(1,4)$.
To avoid the case-by-case analysis as much as possible,
we introduce the following $2$-parameter deformation of $\myAA(2,2)$.
Let $\tilde \FF = \QQ(q_1,q_2,Y_1,Y_2)$ be the field of rational functions
in
four (commuting) independent variables.
We recursively define elements
$Y_m \in \tilde \FF$ for $m \in \ZZ$ by the relations
\begin{equation}
\label{eq:qu-relations}
Y_{m-1} Y_{m+1} = Y_m^2 + q_{\rem{m}} Y_m + 1.
\end{equation}
Now let $\tilde \myAA$ denote the $\ZZ[q_1,q_2]$-subalgebra of $\tilde \FF$
generated by the $Y_m$ for all $m \in \ZZ$.
In view of \cite[Example~5.3]{fz-Laurent},
$\tilde \myAA \subset \ZZ[q_1,q_2,Y_m^{\pm 1}, Y_{m+1}^{\pm 1}]$
for every $m \in \ZZ$, i.e., every element of~$\tilde \myAA$
is a Laurent polynomial in $Y_m$ and $Y_{m+1}$ with coefficients in
$\ZZ[q_1,q_2]$.

The algebra $\myAA(2,2)$ is obtained from $\tilde \myAA$ by an obvious
specialization
\begin{equation}
\label{eq:22-spec}
Y_m  = y_m, \quad q_1 = q_2 = 0.
\end{equation}
Less obvious is the following observation: the subalgebra of
$\myAA (1,4)$ generated by the elements $y_m^{\rem{m}}$ for
$m \in \ZZ$ (i.e., by the $y_m$ for $m$ odd, and the $y_m^2$ for
$m$ even) is obtained from $\tilde \myAA$ by the specialization
\begin{equation}
\label{eq:14-spec}
Y_m = y_m^{\rem{m}}, \quad q_1 = 2, \quad q_2 = 0.
\end{equation}
To see this, notice that
$$Y_{m-1} Y_{m+1} = (y_{m-1} y_{m+1})^2 = (y_m + 1)^2 =
Y_m^2 + 2 Y_m + 1$$
for $m$ odd, and
$$Y_{m-1} Y_{m+1} = y_{m-1} y_{m+1} = y_m^4 + 1 = Y_m^2 + 1$$
for $m$ even.

We introduce an element $Z \in \tilde \myAA$ by setting
\begin{equation}
\label{eq:Z}
Z = Y_0 Y_3 - (Y_1 + q_1)(Y_2 + q_2).
\end{equation}
As before, we define the sequence $Z_1, Z_2, \dots$ of elements of
$\tilde \myAA$ by setting $Z_n = T_n (Z)$. In view of
\eqref{eq:z}, each element $z_n \in \myAA(2,2)$ (resp.~$z_n \in \myAA(1,4)$)
is obtained from $Z_n$ by a specialization \eqref{eq:22-spec}
(resp.~\eqref{eq:14-spec}).

Returning to the cluster algebras of affine types, we now prove the following result, sharpening
Proposition~\ref{pr:z-imaginary} in the spirit of
Proposition~\ref{pr:triangleprop}.
Recall that $\delta$ denotes the minimal positive imaginary root
given by \eqref{eq:delta}.

\begin{proposition}
\label{pr:Newton-zn}
In the affine type, each $z_n$ is a monic Laurent polynomial
in $y_1$ and $y_2$, and its Newton polygon
${\rm Newt}(z_n)$ is equal to the triangle $\Delta(n\delta)$.
\end{proposition}

\begin{proof}
It is enough to show that $Z_n$ is a monic Laurent polynomial
in $Y_1$ and $Y_2$, and that its Newton polygon
${\rm Newt}(Z_n)$ is a triangle $\Delta_n$ with vertices
$(-n,-n), (n,-n)$ and $(-n,n)$; in view of \eqref{eq:delta},
the proposition then follows from the specializations
\eqref{eq:22-spec} and \eqref{eq:14-spec}.

First of all, combining \eqref{eq:Z} with \eqref{eq:qu-relations},
we obtain
\begin{equation}
\label{eq:Z-Laurent}
Z = \frac{Y_1}{Y_2}+\frac{Y_2}{Y_1}+\frac{q_1}{Y_2}
+ \frac{q_2}{Y_1} + \frac{1}{Y_1 Y_2},
\end{equation}
implying our claims for $Z_1 = Z$.
To deal with $n > 1$, we write $Z$ as $Z = t + t^{-1} + u$, where
$$t = \frac{Y_1}{Y_2}, \quad u = \frac{q_1}{Y_2}
+ \frac{q_2}{Y_1} + \frac{1}{Y_1 Y_2} \ .$$
Taking the Taylor expansion of the Chebyshev polynomial $T_n$, we get
\begin{equation}
\label{eq:Tn-Taylor}
Z_n = T_n (Z) = \sum_{k=0}^n \frac{1}{k!} T_n^{(k)}(t+t^{-1}) u^k.
\end{equation}
Since each derivative $T_n^{(k)}$ is a linear combination of the
polynomials $T_\ell$ for $0 \leq \ell \leq n-k$, it follows that
$Z_n$ is a linear combination of Laurent monomials $t^\ell u^k$
for $0 \leq k \leq n$ and $-(n-k) \leq \ell \leq n-k$.
The Newton polygon ${\rm Newt}(t^\ell u^k)$ is a
triangle (possibly degenerate) with vertices
$(\ell - k,-\ell - k), (\ell,-\ell - k)$ and $(\ell - k,-\ell)$,
and sides parallel to those of $\Delta_n$.
By inspection, all these vertices belong to $\Delta_n$, and so
${\rm Newt}(Z_n) \subseteq \Delta_n$.
Furthermore, the only term $t^\ell u^k$ whose Newton polygon may
contain the vertex $(-n,-n)$ (resp.~$(n,-n)$, $(-n,n)$) of
$\Delta_n$ is $u^n$ (resp.~$t^n$, $t^{-n}$), and the corresponding
Laurent monomial occurs with coefficient~1.
But it is clear from \eqref{eq:Tn-Taylor} and the definition of the
Chebyshev polynomial $T_n$ that each of the monomials $u^n$ and $t^{\pm n}$
occurs in
$Z_n$ with coefficient~$1$.
This implies the desired equality ${\rm Newt}(Z_n) = \Delta_n$, as
well as the claim that $Z_n$ is monic.
\end{proof}

A more careful analysis of \eqref{eq:Tn-Taylor} implies the following
analogue of Proposition~\ref{pr:no-positive-terms}.

\begin{proposition}
\label{pr:no-positive-terms-zn}

\begin{enumerate}
\item
Let $n > 0$ and suppose that a point $\gamma \in Q$ belongs to the
interior of the side of the triangle $\Delta(n\delta)$ opposite to the
vertex $- n \delta$.
Then the monomial $y^\gamma$ does not occur in the Laurent expansion
of~$z_n$ in $y_1$ and $y_2$.

\item
The Laurent expansion of $z_n$ in $y_1$ and $y_2$ contains no monomials
$y_1^p y_2^q$
with $p, q \geq 0$.
\end{enumerate}
\end{proposition}

\begin{proof}
As in Proposition~\ref{pr:Newton-zn}, it is enough to prove the
corresponding statements for~$Z_n$.
Using the notation in the proof of Proposition~\ref{pr:Newton-zn},
it is easy to see that the only terms $t^\ell u^k$ whose Newton
polygons have non-empty intersection with the side $[(n,-n), (-n,n)]$ of
$\Delta_n$ are those with $k = 0$, i.e., those appearing in
$T_n (t+t^{-1}) = t^n + t^{-n}$.
This implies Part~1.

Part 2 follows at once from Part~1 by observing that the
intersection of $\Delta_n$ with the positive quadrant $\RR_{\geq 0}^2$
consists of a single point (namely, the origin), and this point
belongs to the interior of the line segment $[(n,-n), (-n,n)]$.
\end{proof}

Our next result implies, in particular, that both Propositions~\ref{pr:Newton-zn} and
\ref{pr:no-positive-terms-zn} hold not just for $\{y_1,y_2\}$ but for any choice of ``initial''
cluster.

\begin{proposition}
\label{pr:z-invariance}
For every $p \in \ZZ$, the automorphism $\sigma_p$ preserves
each of the elements $z_n$.
\end{proposition}

\begin{proof}
It suffices to prove the corresponding statement for the
elements $Z_n \in \tilde \myAA$, and the $\ZZ[q_1, q_2]$-automorphisms
$\sigma_p$ acting on the generators by $\sigma_p (Y_m) = Y_{2p-m}$;
the proposition then follows by an appropriate specialization
(\eqref{eq:22-spec} in the case of $\myAA(2,2)$, and
\eqref{eq:14-spec} in the case of $\myAA(1,4)$).
It is enough to show that $\sigma_1 (Z) = Z$.
This is proved by the following calculation using \eqref{eq:Z-Laurent}:
\begin{align}
\label{eq:Z-invariance}
\sigma_1 (Z) &= \frac{Y_0^2 + q_2 Y_0 + (Y_1^2 + q_1 Y_1 + 1)}{Y_0 Y_1} =
\frac{Y_0^2 + q_2 Y_0 + Y_0 Y_2}{Y_0 Y_1} =
\frac{Y_0 + Y_2 + q_2}{Y_1}\\
\nonumber
&= \frac{Y_0 Y_2+ Y_2^2 + q_2 Y_2}{Y_1 Y_2}
= \frac{Y_1^2 + Y_2^2 + q_1 Y_1 + q_2 Y_2 + 1}{Y_1 Y_2} = Z,
\end{align}
as desired.
\end{proof}

The above proof implies, in particular, that $Z$ can be expressed in any
four consecutive variables $Y_m, \dots, Y_{m+3}$ as
$Z=Y_m Y_{m+3}-(Y_{m+1}+q_{\rem{m+1}})(Y_{m+2}+q_{\rem{m+2}})$ where as above
$\rem{n}$ is $1$ for $n$ odd and $2$ for $n$ even.
In the two affine
cases this specializes for all even $m\in\ZZ$ to
$$
z =
\begin{cases}
\label{eq:explicit-z}
y_m y_{m \pm 3}-y_{m \pm 1}y_{m \pm 2} & \text{if $(b,c) = (2,2)$;} \\[.05in]
y_m^2 y_{m \pm 3} - (y_{m \pm 1} + 2) y_{m\pm 2}^2 & \text{if $(b,c) = (1,4)$.}
\end{cases}
$$

We are ready for the proofs of main results for the affine
types: Theorems~\ref{th:canonicalbasis}--\ref{th:canonicalbasis-symmetry}
and \ref{th:canonicalbasis-affine}.
We follow the same strategy as in Section~\ref{sec:proofs-finite}
with necessary technical modifications.

First of all, Theorem~\ref{th:canonicalbasis-parametrization} for the affine
types is a consequence of Theorem~\ref{th:canonicalbasis-affine}
combined with Propositions~\ref{pr:z-imaginary}
and~\ref{pr:cluster-monomials-leading-terms}.
Second, Propositions~\ref{pr:cluster-monomials-leading-terms}
and \ref{pr:z-imaginary} imply that the cluster monomials together
with the elements $z_n$ are linearly independent: this is proved
by the same argument as in Corollary~\ref{cor:linear-independence}.
Third, Theorem~\ref{th:canonicalbasis-symmetry} follows from
Theorem~\ref{th:canonicalbasis-affine} combined with
Propositions~\ref{pr:sigma-cluster-monomials},
\ref{pr:z-imaginary} and \ref{pr:z-invariance}.

As in Section~\ref{sec:proofs-finite}, it remains to prove the following
counterparts of
\eqref{eq:y-positivity-finite}--\eqref{eq:positive-converse-finite}:
\begin{eqnarray}
\label{eq:yz-positivity-affine}
&&\text{Every cluster variable and every $z_n$ is a positive element of
$\myAA(b,c)$.}\\
\label{eq:spanning-affine}
&&\text{Cluster monomials and the elements $z_n$ linearly span
$\myAA(b,c)$.}\\
\label{eq:positive-converse-affine}
&&\text{If $y \in \myAA(b,c)$ is written as a linear combination of
cluster monomials}\\
\nonumber
&&\text{and the $z_n$'s then its every coefficient is equal to some
coefficient in the}
\\ %\hspace{-.2in}
\nonumber
&&\text{Laurent polynomial expansion of $y$ with respect to some cluster.}
\end{eqnarray}

Our proofs of~\eqref{eq:yz-positivity-affine} and~\eqref{eq:spanning-affine}
will use explicit expressions for all minimal ``forbidden" products
\begin{equation}
\label{eq:forbidden}
z_n z_p \,\, (n, p > 0), \quad z_n y_m \,\, (n > 0, m \in \ZZ),
\quad y_m y_{m+n} \,\, (m \in \ZZ, \ n \geq 2),
\end{equation}
as linear combinations of the cluster monomials and the $z_n$'s.
In the formulas below we use the convention that $z_0 = 1$ and $z_{-n} = 0$
for $n > 0$.
Also $\rem{n}$ has the same meaning as above: it is~$1$ for~$n$ odd, and~$2$
for~$n$ even.

\begin{proposition}
\label{pr:straightening-affine}

\begin{enumerate}
\item
For any affine type, the following relation holds for all $p \geq n \geq 1$:
\begin{equation}
\label{eq:zz}
z_{n} z_{p} =
\begin{cases}
z_{p-n}+z_{p+n} & \text{if $p > n$;} \\[.05in]
2 +z_{2n} & \text{if $p= n$.}
\end{cases}
\end{equation}

\item
In the case $(b,c) = (2,2)$, for all
$ m \in \ZZ$ and $n \geq 1$, we have
\begin{equation}
\label{eq:zy-22}
z_n y_m=y_{m-n}+y_{m+n},
\end{equation}
\begin{equation}
\label{eq:yy-22}
y_{m} y_{m+n} = y_{\lfloor m + \frac{n}{2} \rfloor} y_{\lceil  m +
\frac{n}{2} \rceil}
+ \!\displaystyle \sum_{k \geq 1} k z_{n-2k}.
\end{equation}

\item
In the case $(b,c) = (1,4)$, we have
\begin{equation}
\label{eq:zy-14}
z_n y_m =
\begin{cases}
y_{m-2n} + y_{m+2n} & \text{for $m$ even;}  \\[.05in]
y_{m-n}^{\rem{m-n}} + y_{m+n}^{\rem{m+n}} +
4 \!\displaystyle \sum_{k \geq 1} k z_{n-2k}
& \text{for $m$ odd;}
\end{cases}
\end{equation}
\begin{equation}
\label{eq:yy-14-00}
y_{m} y_{m+2n} = y_{m+n}^{\rem{m+n}} + \!\displaystyle \sum_{k \geq 1}
(2k-1) z_{n+1-2k} \quad \text {($m$ even, $n \geq 0$)};
\end{equation}
\begin{align}
\label{eq:yy-14-01}
y_m y_{m \pm n} =& \!\displaystyle
\sum_{1 < 2k < n} \min(4k,n-2k) \: y_{m \pm 4k} +
\begin{cases}
y_{s}^3 & \text{if $n\equiv 0$  mod $3$,} \\[.05in]
y_{\lfloor s \rfloor} y_{\lceil s \rceil}
& \text{otherwise}
\end{cases}
\\
\nonumber
& \text{($m$ even, $n \geq 1$ odd),}
\end{align}
where $s = m \pm \frac{2n}{3}$;
\begin{align}
\label{eq:yy-14-11}
y_{m} y_{m+2n} =& y_{m+n}^{2\rem{m+n}} +
4 \displaystyle \sum_{k = 1}^{n-1} \min(k,n-k) \: y_{m + 2k}
+ \frac{1}{3} \displaystyle \sum_{k \geq 1} (2k^3 + k) z_{2n-2k}\\
\nonumber
& \text{($m$ odd, $n \geq 0$);}
\end{align}
\end{enumerate}
\end{proposition}

\begin{proof}
To prove Proposition~\ref{pr:straightening-affine}, we first establish
the analogs of the relations in question
for the elements $Z_n$ and $Y_m$ in the algebra $\tilde \myAA$.  As
above, we use the convention that $Z_0 = 1$, and $Z_{-n}=0$ for $n>0$.

\begin{lemma}
\label{lem:ZY-relations}
The following relations hold for all $m \in \ZZ$ and $p \geq n \geq 1$:
\begin{equation}
\label{eq:ZZ}
Z_{n} Z_{p} =
\begin{cases}
Z_{p-n}+Z_{p+n} & \text{if $p > n$;} \\[.05in]
2 +Z_{2n} & \text{if $p= n$;}
\end{cases}
\end{equation}
\begin{equation}
\label{eq:ZY}
Z_n Y_m = Y_{m-n} + Y_{m+n} +
\displaystyle \sum_{k \geq 1} k q_{\rem{m+k}} Z_{n-k};
\end{equation}
\begin{equation}
\label{eq:YY}
Y_{m} Y_{m+n} = Y_{\lfloor m + \frac{n}{2} \rfloor} Y_{\lceil  m +
\frac{n}{2} \rceil}
+ \displaystyle \sum_{k = 1}^{n-1} \min(k,n-k) q_{\rem{m+k}} Y_{m + n- k}
+ \displaystyle \sum_{k \geq 1} c_k Z_{n-1-k},
\end{equation}
where the coefficients $c_k$ are given by
\begin{equation}
\label{eq:YY-coeffs}
c_{2p} = \frac{p(p+1)(2p+1)}{6} q_1 q_2, \quad
c_{2p-1} = \frac{(p-1)p(p+1)}{6} (q_1^2 + q_2^2) + p \ .
\end{equation}
\end{lemma}

\begin{proof}
The relation \eqref{eq:ZZ} follows at once from
the definition of  Chebyshev polynomials $T_n$.

To prove \eqref{eq:ZY}, we start with $n=1$.
Proposition~\ref{pr:z-invariance} and the symmetry of indices $1$ and $2$ in
\eqref{eq:qu-relations} reduces the case of general~$m$
to a special case $m=1$; the corresponding relation
$Z_1 Y_1 = Y_0 + Y_2 + q_2$ has already appeared in \eqref{eq:Z-invariance}.
The case $n =2$ now follows since
$$Z_2 Y_m = (Z_1^2 - 2) Y_m = Z_1 (Y_{m-1} + Y_{m+1} +
q_{\rem{m+1}}) - 2 Y_m = Y_{m-2} + Y_{m+2} + q_{\rem{m+1}} Z_1 +
2 q_{\rem{m}}.$$
Proceeding by induction, we now assume that the products
$Z_{n-1} Y_m$ and $Z_n Y_m$ satisfy \eqref{eq:ZY} for some $n \geq 2$,
and we want to show that the same is true for $Z_{n-1} Y_m$.
This is done by the following calculation using \eqref{eq:ZZ}:
\begin{align*}
&Z_{n+1} Y_m = (Z_1 Z_n - Z_{n-1}) Y_m\\
&= Z_1 (Y_{m-n} + Y_{m+n} +  \displaystyle \sum_{k \geq 1} k q_{\rem{m+k}}
Z_{n-k})
- (Y_{m-n+1} + Y_{m+n-1} +
\displaystyle \sum_{k \geq 1} k q_{\rem{m+k}} Z_{n-1-k})\\
&= Y_{m-n-1} + Y_{m+n+1} + 2 q_{\rem{m+n+1}}
+  \displaystyle \sum_{k \geq 1} k q_{\rem{m+k}} (Z_1 Z_{n-k} -
Z_{n-1-k})\\
&= Y_{m-n-1} + Y_{m+n+1} + (n+1) q_{\rem{m+n+1}}
+  \displaystyle \sum_{k = 1}^n k q_{\rem{m+k}} Z_{n+1-k},
\end{align*}
as desired.

The proof of \eqref{eq:YY} uses the same strategy but requires a
little bit more work.
The case $n=1$ is tautologically true (the two summations in the
right hand side being empty), and the case $n=2$ is the cluster relations
\eqref{eq:qu-relations}.
Proceeding by induction, we now assume that the products
$Y_m Y_{m+n-1}$ and $Y_m Y_{m+n}$ satisfy \eqref{eq:YY} for some $n \geq 2$,
and we want to show that the same is true for $Y_m Y_{m+n+1}$.
Using the relation
$$Y_{m+n+1} = Z_1 Y_{m+n} - Y_{m+n-1} - q_{\rem{m+n+1}},$$
which is a special case of \eqref{eq:ZY}, we obtain
\begin{align*}
&Y_m Y_{m+n+1} = Y_m (Z_1 Y_{m+n} - Y_{m+n-1} - q_{\rem{m+n+1}})\\
&\quad
= Z_1 (Y_{\lfloor m + \frac{n}{2} \rfloor}
Y_{\lceil  m + \frac{n}{2} \rceil}
+ \displaystyle \sum_{k = 1}^{n-1} \min(k,n-k) q_{\rem{m+k}} Y_{m + n- k}
+ \displaystyle \sum_{k \geq 1} c_k Z_{n-1-k})\\
&\qquad
- (Y_{\lfloor m + \frac{n-1}{2} \rfloor} Y_{\lceil  m + \frac{n-1}{2}
\rceil}
+ \displaystyle \sum_{k = 1}^{n-2} \min(k,n-1-k) q_{\rem{m+k}} Y_{m + n-1-
k}
+ \displaystyle \sum_{k \geq 1} c_k Z_{n-2-k})\\
&\qquad - q_{\rem{m+n+1}} Y_m = S_1 + S_2 + S_3,
\end{align*}
where we abbreviated
\begin{align*}
& S_1 = Z_1 Y_{\lfloor m + \frac{n}{2} \rfloor} Y_{\lceil  m + \frac{n}{2}
\rceil}
- Y_{\lfloor m + \frac{n-1}{2} \rfloor} Y_{\lceil  m + \frac{n-1}{2}
\rceil},\\
& S_2 =  \displaystyle \sum_{k = 1}^{n} \min(k,n-k) q_{\rem{m+k}}
(Y_{m + n+1- k} + Y_{m + n-1- k} + q_{\rem{m+n+1+k}}) \\
&\qquad\quad - q_{\rem{m+n+1}} Y_m - \displaystyle \sum_{k = 1}^{n-1}
\min(k,n-1-k) q_{\rem{m+k}} Y_{m + n-1-k},\\
& S_3 = \displaystyle \sum_{k \geq 1} c_k (Z_1 Z_{n-1-k} - Z_{n-2-k}).
\end{align*}
Using \eqref{eq:ZZ}, we simplify
$$S_3 = \displaystyle \sum_{k \geq 1} c_k Z_{n-k} + c_{n-2} - c_n.$$
To simplify  $S_1$ and $S_2$, we first assume that $n = 2p$ is even.
The routine calculations yield
\begin{align*}
&S_1 = Y_{m+p} (Z_1 Y_{m+p} - Y_{m+p-1})  = Y_{m+p} Y_{m+p+1} +
q_{\rem{m+p+1}} Y_{m+p};\\
&S_2 = \displaystyle \sum_{k = 1}^{n} \min(k,n+1-k) q_{\rem{m+k}} Y_{m +
n+1- k}
- q_{\rem{m+p+1}} Y_{m+p} + p^2 q_1 q_2.
\end{align*}
Plugging in the expressions for $c_n$ and $c_{n-2}$ given by the first
equality in
\eqref{eq:YY-coeffs}, we see that
$S_1 + S_2 + S_3$ is indeed the right hand side of \eqref{eq:YY} for $Y_m
Y_{m+n+1}$.
In the case $n = 2p-1$, similar calculations  yield
\begin{align*}
&S_1 = %Y_{m+p-1} (Z_1 Y_{m+p} - Y_{m+p-1})  =
Y_{m+p-1} Y_{m+p+1} +
q_{\rem{m+p+1}} Y_{m+p-1} = Y_{m+p}^2 +  q_{\rem{m+p}} Y_{m+p} +
q_{\rem{m+p+1}} Y_{m+p-1} + 1;\\
&S_2 = \displaystyle \sum_{k = 1}^{n} \min(k,n+1-k) q_{\rem{m+k}} Y_{m +
n+1- k}
-  q_{\rem{m+p}} Y_{m+p}\\
&\qquad \quad- q_{\rem{m+p+1}} Y_{m+p-1} + \frac{p(p-1)}{2}(q_1^2 + q_2^2)
\end{align*}
To conclude that $S_1 + S_2 + S_3$ is again the desired expression, it
remains
to plug in the expressions for $c_n$ and $c_{n-2}$ given by the second
equality in \eqref{eq:YY-coeffs}.
This completes the proof of Lemma~\ref{lem:ZY-relations}.
\end{proof}

Returning to the relations in
Proposition~\ref{pr:straightening-affine},
we notice that all of them in the case $(b,c) = (2,2)$
are obtained from the corresponding relations in
Lemma~\ref{lem:ZY-relations} by a specialization \eqref{eq:22-spec}.
So we assume throughout the rest of the proof that we are in the case
$(b,c) = (1,4)$.

Using the specialization \eqref{eq:14-spec}, we see that
\eqref{eq:zz}, the second case in \eqref{eq:zy-14} and
\eqref{eq:yy-14-11} are again consequences of
the corresponding relations in Lemma~\ref{lem:ZY-relations}.
It remains to prove the first case in \eqref{eq:zy-14}, and also
\eqref{eq:yy-14-00} and \eqref{eq:yy-14-01}.

We start with an observation that the specialization \eqref{eq:14-spec}
of \eqref{eq:ZY} yields the identity
$$z_1 y_{m}^2 = y_{m-1} + y_{m+1} + 2 \quad \text{($m$ even).}$$
Using the exchange relations \eqref{eq:clusterrelations}, we obtain
the following relation for any even $m$:
$$z_1 y_{m} = \frac{y_{m-1} + y_{m+1} + 2}{y_{m}} =
\frac{(y_{m-2}y_{m}-1)  + (y_{m}y_{m+2}-1) + 2}{y_{m}} =
y_{m-2} + y_{m+2}.$$
We see that the operator of multiplication by~$z_1$
preserves the linear span $V$ of the cluster variables~$y_{m}$ with even
$m$;
and it acts on $V$ as $t + t^{-1}$, where $t$ is the shift operator acting
by $t(y_{2m}) = y_{2m+2}$.
It follows that the multiplication by an arbitrary element $z_n$ acts on $V$
as
$T_n (t+ t^{-1}) = t^n + t^{-n}$, which is precisely the first case in
\eqref{eq:zy-14}.

To prove \eqref{eq:yy-14-00}, we proceed by induction on~$n$.
The case $n=0$ is trivial, while $n=1$ is simply the exchange relation.
To treat the product $y_m y_{m + 2(n+1)}$ for $n \geq 1$,
we substitute $y_{m + 2(n+1)} = z_1 y_{m+2n} - y_{m+2(n-1)}$
(see \eqref{eq:zy-14}) and then use the inductive assumption together with
\eqref{eq:zz} and the appropriate specialization of \eqref{eq:ZY};
the routine details are left to the reader.

The same strategy is used for the proof of \eqref{eq:yy-14-01}.
This completes the proof of Proposition~\ref{pr:straightening-affine}.
\end{proof}

Now we are ready to prove~\eqref{eq:yz-positivity-affine},
i.e., that the
cluster variables and the $z_n$ are positive.
By Proposition~\ref{pr:z-invariance}, it suffices
to check the positivity of the Laurent expansion of each $z_n$
or $y_m$ in the variables $y_1$ and $y_2$.
For~$z_n$, the desired positivity is an immediate consequence of
(\ref{eq:Tn-Taylor}) combined
with the following well-known property of Chebyshev polynomials:
the derivative of $T_n$ is a positive linear combination of
the polynomials $T_p$ with $p < n$.
For the sake of completeness, here is the proof:
setting $z = t+t^{-1}$, we obtain
\begin{align*}
d T_n (z)/dz &= d T_n (t + t^{-1})/dt \cdot (dz/dt)^{-1} =
n(t^{n-1} - t^{-n-1})/(1 - t^{-2})\\
&= n(t^n - t^{-n})/(t- t^{-1}) =
n (t^{n-1} + t^{n - 3} + \cdots + t^{-(n-1)})\\
& = n \!\!\!\sum_{0 \leq k \leq \frac{n-1}{2}} \!\!\!\!T_{n-1 -2k}.
\end{align*}

To prove that each cluster variable $y_m$
for $m \geq 3$ is a positive Laurent polynomial in
$y_1$ and $y_2$,  it suffices to apply one of the formulas
(\ref{eq:yy-22}), (\ref{eq:yy-14-00}) or (\ref{eq:yy-14-11})
to the product $y_{\rem{m}} y_m$, and to notice that the
right hand side is a positive linear combination of the
elements~$z_n$ and cluster monomials containing the variables
$y_p$ with $\rem{m} < p < m$.
The desired positivity of $y_m$ follows by induction on~$m$.
The positivity of $y_m$ for $m \leq 0$ follows by
Proposition~\ref{pr:z-invariance}.

\begin{remark}
\label{rem:YZ-positivity}
Note that the same argument establishes that, in the algebra
$\tilde \myAA$ given by \eqref{eq:qu-relations}, every element
$Y_m$ is a Laurent polynomial in $Y_1$ and $Y_2$ whose
coefficients are polynomials in $q_1$ and $q_2$ with \emph{positive}
integer coefficients.
\end{remark}

To prove~\eqref{eq:spanning-affine}, it suffices to show that
every monomial~$M$ in the variables $z_n$ and $y_m$ is a linear
combination of the cluster monomials and the $z_n$'s.
We write~$M$ in the form
$M = z_{n_1}^{a_1} \cdots z_{n_s}^{a_s} y_{m_1}^{b_1} \cdots y_{m_t}^{b_t}$,
where
$0< n_1 < \cdots < n_s, \, m_1 < \cdots < m_t$, and the exponents
$a_1, \dots, a_s$ and $b_1, \dots, b_t$ are positive integers.
We define the \emph{multi-degree}
$\mu(M) = (\mu_1(M), \mu_2(M), \mu_3(M)) \in \ZZ_{\geq 0}^3$ by setting
\begin{align*}
%\label{eq:degree-affine}
&\mu_1(M) =
\begin{cases}
\sum_{i=1}^s a_i + \sum_{j=1}^t b_j
& \text{if $(b,c) = (2,2)$;} \\[.05in]
\sum_{i=1}^s a_i + \sum_{j=1}^t \rem{m_j-1} b_j
& \text{if $(b,c) = (1,4)$;}
\end{cases}\\
& \mu_2 (M) = m_t - m_1;\\
& \mu_3(M) = b_1 + b_t.
\end{align*}
The lexicographic order on $\ZZ_{\geq 0}^3$ makes it into a well
ordered set (i.e., every non-empty subset of $\ZZ_{\geq 0}^3$ has
the smallest element).
Therefore, to complete the proof, it suffices to show that
every monomial $M$ which has at least one of the ``forbidden" products in
\eqref{eq:forbidden}
as a factor, can be written as a linear combination of monomials
of (lexicographically) smaller multi-degree.
We will show that this can be done by replacing
some ``forbidden" factor of $M$ with its expression given by the
appropriate relation in Proposition~\ref{pr:straightening-affine}.
Indeed, if $\sum_{i=1}^s a_i \geq 2$ (resp.~$\sum_{i=1}^s a_i = 1$)
then one can apply \eqref{eq:zz} (resp.~\eqref{eq:zy-22}
or \eqref{eq:zy-14}), expressing $M$ as a linear combination of monomials
with smaller value of $\mu_1$.
So we can assume that $M = y_{m_1}^{b_1} \cdots y_{m_t}^{b_t}$
with $m_t - m_1 \geq 2$.
Then we apply \eqref{eq:yy-22} or one of
\eqref{eq:yy-14-00}--\eqref{eq:yy-14-11}
to the product $y_{m_1} y_{m_t}$.
By inspection, in the resulting expression for~$M$,
there is precisely one monomial $M'$ with
$\mu_1 (M') = \mu_1(M)$, while the rest of the terms have
smaller value of~$\mu_1$.
By further inspection, we have $\mu_2 (M') < \mu_2(M)$
(resp.~$\mu_2 (M') = \mu_2(M)$ and $\mu_3 (M') = \mu_3 (M)-2$)
if $\min(b_1,b_t) = 1$ (resp.~$\min(b_1,b_t) \geq 2$).
This concludes the proof of \eqref{eq:spanning-affine}.

To prove \eqref{eq:positive-converse-affine}, suppose that an element
$y \in \myAA(b,c)$ is
written (uniquely) as a linear combination of cluster monomials and the
elements $z_n$.
We need to show that every coefficient in this linear combination
is equal to some coefficient in the Laurent expansion of $y$ with
respect to some cluster.
For the coefficients of cluster monomials,
this is proved by precisely the same argument as in the proof of
\eqref{eq:positive-converse-finite} given in
Section~\ref{sec:proofs-finite}; the only extra thing to take into
account is Part~2 of Proposition~\ref{pr:no-positive-terms-zn}.
It remains to deal with the coefficient of each $z_n$.
By Proposition~\ref{pr:z-invariance}, we can assume without loss
of generality that our linear combination
involves only cluster variables $y_m$ with $m \geq 3$.
It is then enough to show the following: for every $n \geq 1$,
there exists a Laurent monomial $y^\gamma = y_1^{g_1} y_2^{g_2}$ that occurs
with coefficient~$1$ in (the Laurent expansion of) $z_n$ but does
not occur in any $z_p$ for $p \neq n$, or in any cluster monomial
in the variables $y_m$ for $m \geq 3$.
We claim that the following vector $\gamma$ has all the desired properties:
$$\gamma = (g_1,g_2) =
\begin{cases}
(n, - n) & \text{if $(b,c) = (2,2)$;} \\[.05in]
(n, - 2n) & \text{if $(b,c) = (1,4)$.}
\end{cases}$$
Indeed, by Proposition~\ref{pr:Newton-zn}, $\gamma$ is a vertex of
${\rm Newt}(z_n)$, and so $y^\gamma$ occurs in $z_n$ with
coefficient~$1$.
The fact that $y^\gamma$ does not occur in $z_p$ for $p \neq n$
follows from Part~1 of Proposition~\ref{pr:no-positive-terms-zn}.
Finally, the fact that $y^\gamma$ does not occur in any cluster
monomial $y_m^p y_{m+1}^q$ for $m \geq 3$, follows
from Corollary~\ref{cor:phi-affine}: an immediate check shows that
$\gamma$ does \emph{not} satisfy the linear inequality given there
that must hold on ${\rm Newt}(y_m^p y_{m+1}^q)$.
This concludes the proofs of Theorems~\ref{th:canonicalbasis} and
\ref{th:canonicalbasis-affine} for the affine type.

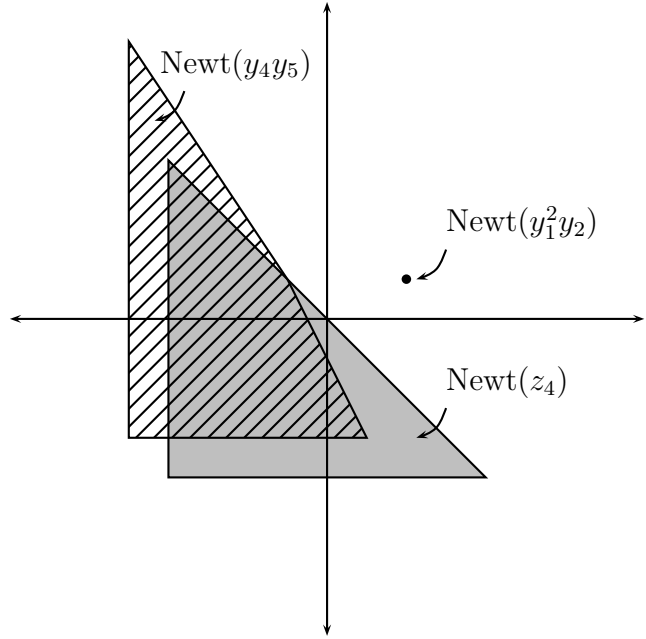
\begin{figure}[ht]
\begin{center}
\psset{unit=15pt}
\psset{labelsep=4pt}
\newgray{lightgray}{0.75}
\begin{pspicture}(-8,-8)(8,8)
\pspolygon[fillcolor=lightgray,fillstyle=solid](-4,4)(4,-4)(-4,-4)
\pspolygon[fillstyle=hlines](-5,-3)(1,-3)(-1,1)(-5,7)
\psdots[dotscale=1.0](2,1)
\psline[arrows=<->](-8,0)(8,0)
\psline[arrows=<->](0,-8)(0,8)
\pscurve[arrows=->](3,-2.25)(2.75,-2.75)(2.25,-3)
\pscurve[arrows=->](-3.6,5.75)(-3.85,5.25)(-4.35,5)
\pscurve[arrows=->](3,1.75)(2.75,1.25)(2.25,1)

\uput[90](3,-2.25){\put(0,0){Newt($z_4$)}}
\uput[90](-3.85,5.25){\put(-.35,.7){Newt($y_4y_5$)}}
\uput[90](3,1.75){\put(0,0){Newt($y_1^2y_2$)}}

\end{pspicture}
\end{center}
\caption{Newton polygons of $y_1^2y_2, ~y_4y_5$ and $z_4$
  as Laurent polynomials in $\{y_1,y_2\}$ for $b=c=2$.}
\label{fig:posexample-y}
\end{figure}

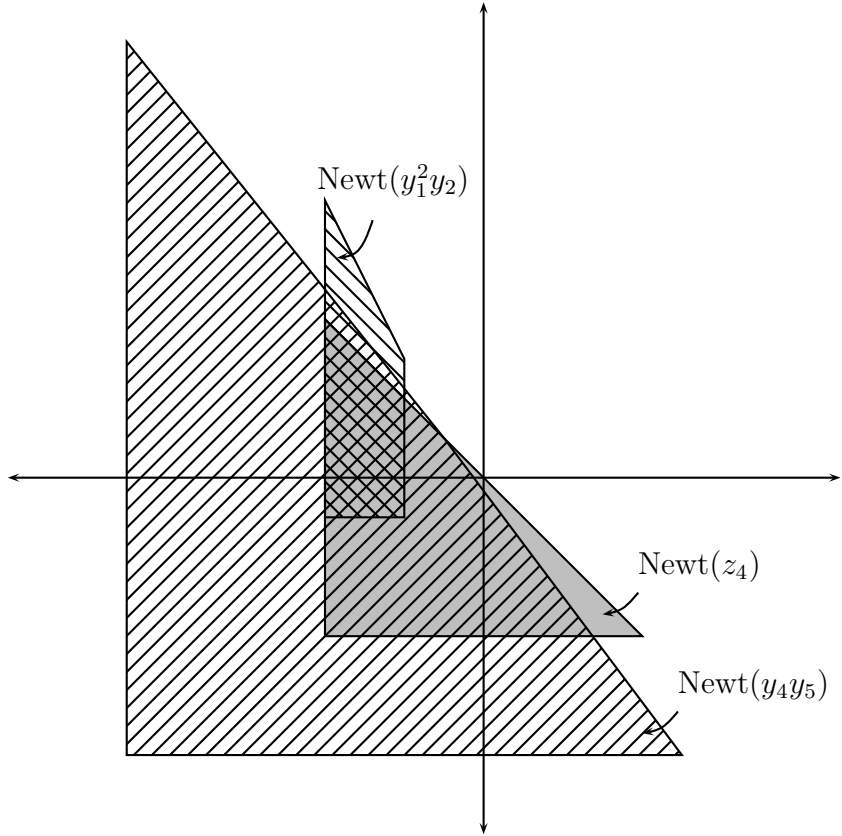
\begin{figure}[ht]
\begin{center}
\psset{unit=15pt}
\psset{labelsep=4pt}
\newgray{lightgray}{0.75}
\begin{pspicture}(-12,-9)(9,12)
\pspolygon[fillcolor=lightgray,fillstyle=solid](-4,4)(4,-4)(-4,-4)
\pspolygon[fillstyle=hlines](-9,-7)(-9,11)(-1,1)(5,-7)
\pspolygon[fillstyle=vlines](-4,-1)(-2,-1)(-2,3)(-4,7)
\psline[arrows=<->](-12,0)(9,0)
\psline[arrows=<->](0,-9)(0,12)
\pscurve[arrows=->](3.9,-2.9)(3.5,-3.3)(3,-3.45)
\uput[90](3.9,-2.9){\put(0,0){Newt($z_4$)}}
\pscurve[arrows=->](-2.8,6.5)(-3.2,5.7)(-3.7,5.55)
\uput[90](-3.0,6.5){\put(-1.2,.5){Newt($y_1^2y_2$)}}
\pscurve[arrows=->](4.9,-5.9)(4.5,-6.3)(4,-6.45)
\uput[90](4.9,-5.9){\put(0,0){Newt($y_4y_5$)}}
\end{pspicture}
\end{center}
\caption{Newton polygons of $y_1^2y_2,~ y_4y_5$ and $z_4$
  as Laurent polynomials in $\{y_{-1},y_0\}$ for $b=c=2$.}
\label{fig:posexample-z}
\end{figure}

\begin{example}
As an illustration of the above proof of \eqref{eq:positive-converse-affine},
let $(b,c)=(2,2)$ and let $y$ be a
positive element that is a $\ZZ$-linear combination of $y_1^2y_2$,
$y_4y_5$ and $z_4$. We wish to show each of the three coefficients is
positive.  First we consider $y$ as a Laurent polynomial in the
variables $\{y_1,y_2\}$.
Propositions~\ref{pr:triangleprop} and \ref{pr:Newton-zn} make it an easy exercise to
compute the Newton polygons of the three terms; the results are shown
in Figure~\ref{fig:posexample-y}.  Since the point ${\rm Newt} (y_1^2y_2)$
does not belong to ${\rm  Newt} (y_4y_5)$ or ${\rm  Newt} (z_4)$, the coefficient of $y_1^2y_2$ must be positive.
Expressing $y$ as Laurent polynomial in $\{y_4,y_5\}$ we could generate
a similar picture (not shown) implying $y_4y_5$ has positive
coefficient.  Finally, we consider $y$ as a Laurent polynomial in
$\{y_{-1},y_0\}$ as shown in Figure \ref{fig:posexample-z}.  Since
$z_4$ is monic and ${\rm Newt}(z_4)$ has a vertex that lies outside
${\rm Newt}(y_1^2y_2)$ and ${\rm Newt}(y_4y_5)$, we conclude that $z_4$ has a
positive coefficient.
\end{example}

\begin{remark}
\label{rem:positivity-non-local}
The second equality in \eqref{eq:upper} shows that the
Laurentness property is \emph{local}: if an element of the ambient
field $\FF$ of the cluster algebra $\myAA (b,c)$
is a Laurent polynomial with respect to three consecutive
clusters, then the same is true with respect to \emph{any} cluster.
Unfortunately, there seems to be no hope to extend this to the
``positive Laurentness" property, as shown by the following
counterexample for the algebra $\myAA(2,2)$.
A direct check shows that the element $y_0 y_1 + y_2 y_3 + y_3 y_4
- z$ has positive Laurent expansion in $y_1$ and $y_2$.
It follows that, for every $n \geq 1$, the element
$$\sum_{m=0}^{n+2} y_m y_{m+1} - z$$
has positive Laurent expansion in $y_m$ and $y_{m+1}$, for
each $m = 1, 2, \dots, n$.
However, by Theorem~\ref{th:canonicalbasis}, this element is not
positive.
\end{remark}

\section{General coefficients}
\label{sec:coeffs-removed}

According to \cite{fz-ClusterI}, for a given pair $(b,c)$ of
positive integers, the most general form of the corresponding
algebra involves the exchange relations
\begin{equation}
\label{eq:clusterrelations-coeffs}
\hat y_{m-1} \hat y_{m+1} =
\left\{
\begin{array}[h]{ll}
q_m \hat y_m^b + r_m & \quad \mbox{$m$ odd;} \\
q_m \hat y_m^c + r_m & \quad \mbox{$m$ even.}
\end{array}
\right.
\end{equation}
Here the coefficients $q_m$ and $r_m$ are the generators of the
\emph{universal coefficient group}~$\PP$ defined by the relations
\begin{equation}
\label{eq:coeffs-relations}
r_{m-1} r_{m+1} =
\left\{
\begin{array}[h]{ll}
q_{m-1} q_{m+1} r_m^c & \quad \mbox{$m$ odd;} \\
q_{m-1} q_{m+1} r_m^b  & \quad \mbox{$m$ even.}
\end{array}
\right.
\end{equation}
The group $\PP$ is a free abelian group with countably many
generators; as a set of free generators, one can choose
$\{q_m : m \in \ZZ\} \cup \{r_1\}$
(cf.~\cite[Proposition~5.2]{fz-ClusterI}).
For our current purposes, we prefer to work with the completion
$\hat \PP$ of $\PP$ given by
\begin{equation}
\label{eq:Phat}
\hat \PP = \PP \otimes_\ZZ \QQ \ ;
\end{equation}
in plain words, $\hat \PP$ is obtained from $\PP$ by adjoining the
roots of all degrees from all the elements of $\PP$.
Thus, we will work over the ground ring $\ZZ \hat \PP$, the
integer group ring of $\hat \PP$.
And so we define the \emph{completed universal cluster algebra}
$\hat \myAA (b,c)$ as the $\ZZ \hat \PP$-subalgebra generated by
the $y_m$ for $m \in \ZZ$ inside the ambient field of rational
functions in $\hat y_1$ and $\hat y_2$ over the field of fractions
of~$\ZZ \hat \PP$.
The \emph{coefficient-free} cluster algebra $\myAA(b,c)$ studied
above is obtained from $\hat \myAA (b,c)$ by the specialization
\begin{equation}
\label{eq:coeff-free-spec}
\hat y_m = y_m, \quad q_m = r_m = 1.
\end{equation}
The following result shows that this specialization does not lead
to any loss of information, since $\hat \myAA (b,c)$ can be
obtained from $\myAA (b,c)$ by an extension of scalars from $\ZZ$
to $\ZZ \hat \PP$.

\begin{proposition}
\label{pr:coeffs-removed}
There is a natural $\ZZ \hat \PP$-algebra isomorphism
$\psi: \myAA (b,c) \otimes \ZZ \hat \PP \to \hat \myAA (b,c)$
given by
\begin{equation}
\label{eq:rescaling}
\psi (y_m) =
\left\{
\begin{array}[h]{ll}
\left (\frac{q_{m}}{r_m}\right)^{1/b} \hat y_m & \quad \mbox{$m$ odd;} \\
\left (\frac{q_{m}}{r_m}\right)^{1/c} \hat y_m & \quad \mbox{$m$ even.}
\end{array}
\right.
\end{equation}
\end{proposition}

\begin{proof}
We just need to show that the elements in \eqref{eq:rescaling}
satisfy the coefficient-free exchange relations
\eqref{eq:clusterrelations}.
This follows at once from \eqref{eq:clusterrelations-coeffs}
and \eqref{eq:coeffs-relations}.
\end{proof}

In view of Proposition~\ref{pr:coeffs-removed}, the canonical
basis $\BB$ in $\myAA (b,c)$ given by Theorem~\ref{th:canonicalbasis}
can be viewed as a $\ZZ \hat \PP$-basis in $\hat \myAA (b,c)$.
Theorem~\ref{th:canonicalbasis} extends to $\hat \myAA (b,c)$ in
the following way.
Let $\ZZ_{> 0} \hat \PP$ denote the semiring in $\ZZ \hat \PP$
consisting of positive integer linear combinations of elements of $\hat \PP$.
The following definition is a counterpart of
Definition~\ref{def:positivity}.

\begin{definition}
\label{def:positivity-coeff}
A non-zero element $y \in \hat \myAA(b,c)$ is \emph{positive} if
for every $m \in \ZZ$, all the coefficients
in the expansion of $y$ as a Laurent polynomial in $\hat y_m$ and $\hat y_{m+1}$
belong to $\ZZ_{> 0} \hat \PP$.
\end{definition}

A counterpart of Theorem~\ref{th:canonicalbasis} can be now stated as follows.

\begin{theorem}
\label{th:canonicalbasis-coeff}
Suppose that $bc \leq 4$.
Then the semiring of positive elements of $\hat \myAA(b,c)$ consists
precisely of $\ZZ_{> 0} \hat \PP$-linear combinations of elements of~$\BB$.
This property determines $\BB$ uniquely up to rescaling by elements of~$\hat \PP$.
\end{theorem}

\begin{proof}
For the finite (resp.~affine) type, the theorem follows at once from the properties
\eqref{eq:y-positivity-finite} and \eqref{eq:positive-converse-finite}
(resp.~\eqref{eq:yz-positivity-affine} and
\eqref{eq:positive-converse-affine}).
\end{proof}

\section*{Acknowledgments}
We thank Sergey Fomin for editorial suggestions, and
Jim Propp and Dylan Thurston for helpful discussions.
Dylan's remark that Chebyshev polynomials of the first kind
appeared naturally in his topological model of cluster algebras was
especially stimulating.
One of the authors (A.~Z.) gratefully acknowledges the hospitality of David
Eisenbud and the financial support of MSRI during his visit in May~2003.

Andrei Zelevinsky is especially happy to dedicate this paper to
Borya Feigin, an old friend and classmate since 7th grade in the
unforgettable Moscow School~No.~2.
Many years ago, Borya taught me the basics of Kac-Moody algebras and their
root systems (he did it in the best possible way, by engaging me in
writing a joint survey paper~\cite{feiginz} on the subject).
I am happy to be able to apply some of his lessons in this work.

\end{document}